\documentclass[12pt]{amsart}

\usepackage{amsmath,amssymb,amsthm,amscd,a4wide,upref}
\usepackage{mathrsfs}
\usepackage{graphicx}
\usepackage[colorlinks=true, linkcolor=blue, citecolor=magenta, urlcolor=blue, backref=page]{hyperref}
\hypersetup{colorlinks,citecolor=blue,filecolor=black,linkcolor=blue,urlcolor=blue}
\usepackage{enumerate}
\usepackage{mathtools}
\usepackage{enumitem}
\usepackage{lmodern}
\usepackage[T1]{fontenc}

\newcommand{\la}{\lambda}
\newcommand{\La}{\Lambda}

\newcommand{\ga}{\gamma}

\newcommand{\de}{\delta}
\newcommand{\ze}{\zeta}

\newcommand{\ts}{\,}

\newcommand{\U}{ {\rm U}}

\newcommand{\Y}{ {\rm Y}}

\newcommand{\gl}{\mathfrak{gl}}

\newcommand{\h}{\mathfrak h}

\newcommand{\sgn}{ {\rm sgn}\ts}

\newcommand{\bi}{\bar{\imath}}
\newcommand{\bj}{\bar{\jmath}}
\newcommand{\bk}{\bar{k}}
\newcommand{\bl}{\bar{l}}

\newcommand{\xiL}{\xi_{\Lambda}}
\newcommand{\chil}{\chi_{i,\La}}
\newcommand{\chjl}{\chi_{m+j,\La}}
\newcommand{\chkll}{\chi_{k,\Lambda^{(1)},\ldots,\Lambda^{(N)}}}

\newcommand{\zeil}{\zeta_{i,\Lambda}}
\newcommand{\zejl}{\zeta_{m+j,\Lambda}}

\newcommand{\zekll}{\zeta_{k,\Lambda^{(1)},\ldots,\Lambda^{(N)}}}

\newcommand{\hs}{h^{(s)}}
\newcommand{\xiLL}{\xi_{\Lambda^{(1)},\ldots,\Lambda^{(N)}}}

\newtheorem{thm}{Theorem}[section]
\newtheorem{lem}[thm]{Lemma}
\newtheorem{prop}[thm]{Proposition}
\newtheorem{cor}[thm]{Corollary}
\newtheorem{conj}[thm]{Conjecture}

\theoremstyle{definition}
\newtheorem{defin}[thm]{Definition}

\theoremstyle{remark}
\newtheorem{remark}[thm]{Remark}
\newtheorem{example}[thm]{Example}

\newcommand{\bth}{\begin{thm}}
\renewcommand{\eth}{\end{thm}}
\newcommand{\bpr}{\begin{prop}}
\newcommand{\epr}{\end{prop}}
\newcommand{\ble}{\begin{lem}}
\newcommand{\ele}{\end{lem}}
\newcommand{\bco}{\begin{cor}}
\newcommand{\eco}{\end{cor}}
\newcommand{\bde}{\begin{defin}}
\newcommand{\ede}{\end{defin}}
\newcommand{\bex}{\begin{example}}
\newcommand{\eex}{\end{example}}
\newcommand{\bre}{\begin{remark}}
\newcommand{\ere}{\end{remark}}
\newcommand{\bcj}{\begin{conj}}
\newcommand{\ecj}{\end{conj}}

\newcommand{\bal}{\begin{aligned}}
\newcommand{\eal}{\end{aligned}}
\newcommand{\beq}{\begin{equation}}
\newcommand{\eeq}{\end{equation}}
\newcommand{\ben}{\begin{equation*}}
\newcommand{\een}{\end{equation*}}

\newcommand{\bpf}{\begin{proof}}
\newcommand{\epf}{\end{proof}}
\newcommand{\Ymn}{\Y(\mathfrak{gl}_{m|n})}            

\newcommand{\BC}{\mathbb{C}}            
\newcommand{\BZ}{\mathbb{Z}}            

\allowdisplaybreaks

\def\beql#1{\begin{equation}\label{#1}}

\begin{document}

\title[Representations of Super Yangians with Gelfand-Tsetlin bases]{Representations of Super Yangians with Gelfand-Tsetlin bases}
\author{Vyacheslav Futorny}
\address{Shenzhen International Center for Mathematics, Southern University of Science and Technology, Shenzhen, China}
\email{vfutorny@gmail.com}
\author{Zheng Li}
\address{School of Artificial Intelligence, Jianghan University, Wuhan, Hubei 430056, China}
\email{lz1994@jhun.edu.cn}
\author{Jian Zhang}
\address{School of Mathematics and Statistics,  Central China Normal University, Wuhan, Hubei 430079, China}
\email{jzhang@ccnu.edu.cn}
\thanks{{\scriptsize
\hskip -0.6 true cm MSC (2020): Primary: 17B37; Secondary: 81R10.
}}
\begin{abstract}
The evaluation homomorphism from the super Yangian $\Ymn$ to the universal enveloping algebra $\U(\gl_{m|n})$ allows one to regard the covariant tensor modules of $\gl_{m|n}$ as $\Ymn$ modules.
In this paper, we derive necessary and sufficient conditions for a simple tensor product of covariant evaluation $\Ymn$ modules to be tame.
We also give an explicit construction of Gelfand-Tsetlin basis for tensor products of skew modules.

\end{abstract}

\maketitle

\section{Introduction}
Consider the following chain of subalgebras of Yangian $\Y(\mathfrak{gl}_m)$ ,
\[
\Y(\mathfrak{gl}_1)\subset\Y(\mathfrak{gl}_2)\subset\cdots\subset\Y(\mathfrak{gl}_m).
\]
The Gelfand-Tsetlin subalgebra  $A(\mathfrak{gl}_m)$ of  $\Y(\mathfrak{gl}_m)$ is generated by the centers of subalgebras in the chain, which is maximal commutative (see \cite{Cher}, \cite{NT1}).
The generators of this subalgebra can be described by quantum determinants \cite{MNO,Mol1}.
A finite-dimensional module with a semisimple action of the subalgebra $A(\mathfrak{gl}_m)$ is called tame.
A finite dimensional irreducible module $V$ is tame if and only if the module $V$ is obtained by pulling back through some automorphism $\omega_f$ from the tensor product
\[
V_{z_1}(\la_1\slash\mu_1)\otimes V_{z_2}(\la_2\slash\mu_2)\otimes \cdots\otimes V_{z_k}(\la_k\slash\mu_k)
\]
of the elementary \( \Y(\mathfrak{gl}_m) \)-modules, for some skew Young diagrams \( \la_1\slash\mu_1, \ldots, \la_k\slash\mu_k \) and some complex numbers \( z_1, \ldots, z_k \) such that \( z_i - z_j \notin \mathbb{Z} \) when \( i \neq j \). This  was
conjectured by I. Cherednik and proved in \cite{Cher} under certain extra conditions
on the module \( V \). In full generality this fact was proved in \cite{NT2}.
It turns out that the spectrum of the action of the subalgebra \( A(\mathfrak{gl}_n) \) in every irreducible tame module is simple.
The eigenbases of \( A(\mathfrak{gl}_n) \) in the tame modules are called Gelfand-Tsetlin bases. Such bases were described in \cite{Mol,NT1,NT2}.

Super Yangian $\Ymn$ was introduced by Nazarov in \cite{Naz1}.
Its simple finite-dimensional modules were characterized by Drinfeld polynomials \cite{Zhang1,Zhang2}.
Consider the following chain of subalgebras of Yangian $\Ymn$,
\[\Y(\mathfrak{gl}_1)  \subset \cdots \subset \Y(\mathfrak{gl}_m)\subset \Y(\mathfrak{gl}_{m|1})\subset \Y(\mathfrak{gl}_{m|n}).\]
The  Gelfand-Tsetlin subalgebra  \( A(\gl_{m|n}) \) of
$\Ymn$ is generated by the centers of subalgebras in the chain. A finite-dimensional module with a semisimple action of the subalgebra  \( A(\gl_{m|n}) \) is called tame.

Similar to $\Y(\gl_n)$, the evaluation homomorphism from the super Yangian $\Y(\gl_{m|n})$ to the universal enveloping algebra $\U(\gl_{m|n})$ allows one to view $\gl_{m|n}$-modules as $\Y(\gl_{m|n})$-modules.
In contrast to the case of $\Y(\gl_n)$, not every finite-dimensional simple module $L(\lambda)$ of $\gl_{m|n}$ is tame.
There are two known types of tame modules: covariant tensor modules and essentially typical modules; see \cite{SV,Mol2,Palev1989a,Palev1989b}.
A covariant tensor module $L(\lambda)$ of $\gl_{m|n}$ can be realized as a submodule of a tensor product of natural representations; see \cite{BR, Ser}.
The Gelfand-Tsetlin bases for covariant tensor modules were described in \cite{SV} and \cite{Mol2}.
When regarded as modules over the super Yangian $\Y(\gl_{m|n})$, covariant tensor modules will be called covariant evaluation modules.

Let $\la^{(i)}$ be a covariant $\gl_{m+r|n}$-weight, $\mu^{(i)}$ be a $\gl_r$-weight, $r$ be a non-negative integer, $N$ be a positive interger and $h^{(i)}\in\mathbb{C}$ for $1\leq i\leq N$.
Let us consider the product of skew modules
\begin{equation}\label{tensor product skew modules}
W=L_{h^{(1)}}(\la^{(1)}\slash\mu^{(1)})\otimes\cdots\otimes L_{h^{(N)}}(\la^{(N)}\slash\mu^{(N)}).
\end{equation}
Via the Hopf algebra structure of $\Ymn$ together with the evaluation homomorphism,
\eqref{tensor product skew modules} can be regarded as a $\Ymn$-module.
In \cite{LM}, Lu and Mukhin proved that \eqref{tensor product skew modules} is a simple $\Ymn$-module if $h^{(s)}-h^{(t)}\notin\BZ$ by using the Drinfeld functor.
In \cite{Lu}, Lu constructed a Gelfand-Tsetlin type basis for $L_h(\la\slash \mu)$ and proved that $W$ is tame and simple in the special case $N=1$.

In this paper, we construct a Gelfand-Tsetlin type basis of \eqref{tensor product skew modules} for any positive integer $N$ satisfying the condition $h^{(s)}-h^{(t)}\notin\BZ$. As a consequence, we show that \eqref{tensor product skew modules} is simple and tame under this condition.
This result serves as a motivation for our subsequent study of tensor products of covariant evaluation modules.

Consider the simple tensor product of covariant evaluation modules, i.e.,
\[
V=L_{h^{(1)}}(\la^{(1)})\otimes L_{h^{(2)}}(\la^{(2)})\otimes\cdots\otimes L_{h^{(M)}}(\la^{(M)}),
\]
such that every $\la^{(i)}$ ($1\leq a\leq M$) is a covariant $\gl_{m|n}$-weight.
In \cite{FLZ}, we provided the irreducibility criterion for tensor products of covariant evaluation modules. In this paper,
we study the necessary and sufficient conditions for simple  tensor products of covariant evaluation modules to be tame. 
We prove that a simple tensor product of covariant evaluation modules is tame if and only if the tensor product of every pair of its factors is tame (Theorem \ref{thm:tame_multiple}).
For the simple tensor product of two covariant evaluation modules,
we give an explicit tameness criterion  in terms of their highest weights and evaluation parameters (Theorem\ref{thm:ns con}).
Furthermore, we establish the connection between the simple tensor product
\[
V_{h}(\underline{\la})=L_{h}(\la^{(1)})\otimes L_{h}(\la^{(2)})\otimes\cdots\otimes L_{h}(\la^{(M)})
\]
and the skew module.
When $n=0$, $V_{h}(\underline{\la})$ is indeed the skew module from \cite{NT2}.

This paper is organized as follows, in Section 2, we review some known facts for super Yangian $\Ymn$.
We omit the proofs of some well-known results, while providing proofs for those results that may not be widely known.
In Section 3, we recall the Gelfand-Tsetlin basis of a covariant evaluation $\Ymn$-module $L(\la)$ and the skew module of $\Ymn$.
Then we give a Gelfand-Tsetlin type basis of the tensor product of skew modules in Section 4.
It implies that such modules are tame and simple.
In Section 5, we study the necessary and sufficient condition for simple tensor products of covariant evaluation $\Ymn$-modules to be tame.
In Section 6, a connection between the module constructed in Section 5 and the tensor product of skew modules in Section 4 is established.

\vspace{1em}

\section{Super Yangians}
\subsection{Lie superalgebra}
Throughout the paper, we work on $\BC$. A super vector space $\BC^{m|n}$ is a ${\BZ}_2$-graded vector space.
Vectors in the $\bar{0}$-graded part $({\BC}^{m|n})_{\bar{0}}={\BC}^m=\sum_{i=1}^m {\BC}
v_i$ are called even, while those in the $\bar{1}$-graded part $({\BC}^{m|n})_{\bar{1}}={\BC}^n=\sum_{j=m+1}^{m+n}{\BC} v_j$ are called odd.
We define the parity $\bi$ of $i$ by
\beql{parity}
\bi =
\begin{cases}
0, & \text{if } i \leq m, \\
1, & \text{if } i > m.
\end{cases}
\eeq

The Lie superalgebra $\gl_{m|n}$ is generated by elements $E_{ij}$, $1\leq i,j\leq m+n$, with the supercommutator relations
\[
[E_{ij},E_{kl}]=\delta_{jk}E_{il}-(-1)^{(\bi+\bj)(\bk+\bl)}\delta_{il}E_{kj},
\]
where the parity of $E_{ij}$ is $\bi+\bj$. Set $e_{i}:=E_{i,i+1}$ and $f_i:=E_{i+1,i}$ for $1\leq i\leq m+n-1$. Denote by $\U(\gl_{m|n})$ the universal enveloping superalgebra of $\gl_{m|n}$.

The Cartan subalgebra $\h$ of $\gl_{m|n}$ is spanned by $E_{ii}$, $1\leq i\leq m+n$. Let $\epsilon_i$, $1\leq i\leq m+n$, be a basis of $\h^*$ (the dual space of $\h$) such that $\epsilon_i(E_{jj})=\delta_{ij}$. There is a bilinear form $(\ ,\ )$ on $\h^*$ given by $(\epsilon_i,\epsilon_j)=(-1)^{\bi}\delta_{ij}$. Define the simple roots $\alpha_i:=\epsilon_i-\epsilon_{i+1}$, for $1\leq i\leq m+n-1$.

Let $\la=(\la_1,\la_2,\dots,\la_{m+n})$ be a tuple of complex numbers. We call $\la$ a $\gl_{m|n}$-weight. Denote $L(\la)$ the irreducible module of $\gl_{m|n}$ generated by a nonzero vector $v$ satisfying the conditions
\[
E_{ii}v=\la_i v,\quad E_{jk}v=0,
\]
for $1\leq i\leq m+n$ and $1\leq j< k\leq m+n$.

It is clear that $\gl_{m|n}$ has a natural representation on $\BC^{m|n}$ such that $E_{ij}v_k=\delta_{jk}v_i$, called vector representation. The highest weight of $\BC^{m|n}$ is the tuple $(1,0,\dots,0)$.

We call $\la=(\la_1,\la_2,\ldots,\la_{m+n})$ a covariant $\gl_{m|n}$-weight,
if $\la$ satisfies:
\begin{enumerate}
\item $\la_1,\ldots,\la_{m+n}$ are all nonnegative integers,
\item $\la_1\geq \ldots\geq \la_m$ and $\la_{m+1}\geq\ldots\geq\la_{m+n}$,
\item $\sharp\{\la_{m+j}>0|1\leq j\leq n\}\leq \la_m$.
\end{enumerate}

We call $L(\la)$ a \emph{covariant module} if $\la$ is a covariant $\gl_{m|n}$-weight. Note that in this case $L(\la)$ is a submodule of $(\BC^{m|n})^{\otimes |\la|}$, where $|\la|=\sum_{i=1}^{m+n}\la_i$.

\subsection{Definition of super Yangian}\label{sec:def}
The super Yangian $\Ymn$ associated with $\gl_{m|n}$ is a $\mathbb{Z}_2$-graded associative algebra over $\mathbb{C}$ with generators $t_{ij}^{(r)}$ for $1 \leq i,j \leq m+n$ and $r \geq 1$.
The defining relations are
\beql{def-gen}
[t_{ij}^{(r)},t_{kl}^{(s)}]=(-1)^{\bi\bj+\bi\bk+\bj\bk}\sum_{p=0}^{\text{min}\{r,s\}-1}(t_{kj}^{(p)}t_{il}^{(r+s-1-p)}-t_{kj}^{(r+s-1-p)}t_{il}^{(p)}),
\eeq
where the parity $\bi$ of $i$ is defined by \eqref{parity}.

Define the formal power series
\[
t_{ij}(u) = \delta_{ij} + t_{ij}^{(1)}u^{-1} + t_{ij}^{(2)}u^{-2} + \cdots.
\]
The defining relations of $\Ymn$ can be written as
\beql{def-rel}
(u-v)[t_{ij}(u),t_{kl}(v)]=(-1)^{\bi\bj+\bi\bk+\bj\bk}(t_{kj}(u)t_{il}(v)-t_{kj}(v)t_{il}(u)).
\eeq
Here we write bracket for the super-commutator.

Let $T(u) = \sum_{i,j=1}^{m+n} t_{ij}(u) \otimes E_{ij}(-1)^{\bj(\bi+1)}$.
Then \eqref{def-rel} can be expressed as the RTT relation form:
\beql{RTT}
R(u-v) T_1(u) T_2(v) = T_2(v) T_1(u) R(u-v),
\eeq
where $R(u) = 1 - \frac{P}{u}$ is the $R$-matrix and $P=\sum_{i,j=1}^{m+n}E_{ij}\otimes E_{ji}(-1)^{\bj}$ is the graded permutation operator.

$\Ymn$ is a Hopf algebra with coproduct
\beql{copro}
\Delta : t_{ij}(u)  \mapsto \sum_{k=1}^{m+n} t_{ik}(u)  \otimes t_{kj}(u),
\eeq
the counit and antipode are defined by
\[
\epsilon: T(u)\rightarrow 1 \text{ \ and \ } S(T(u))=T^{-1}(u),
\]
respectively.

%

\subsection{Berezinian}
Throughout this paper, we use the following notation for  the inverse of $T(u)$
\[
T(u)^{-1}=\left(t'_{ij}(u)\right)_{i,j=1}^{m+n}.
\]


Define the parameter $\ga_k$ by
\[
\ga_k =
\begin{cases}
-k+1, & \text{if } k \leq m, \\
-2m+k, & \text{if } k > m.
\end{cases}
\]
Denote the $k$-th symmetric group by $\mathcal{S}_k$.
\bde
The Berezinian of $\Ymn$ is defined by
\begin{equation}\label{bere}
\begin{split}
B(u) = &\sum_{\sigma \in S_m} \sgn(\sigma) \, t_{\sigma(1)1}(u+\ga_1)t_{\sigma(2)2}(u+\ga_2)  \cdots t_{\sigma(m)m}(u+\ga_m)\\
&\times\sum_{\tau \in S_n} \sgn(\tau) \, t'_{m+1,m+\tau(1)}(u+\ga_{m+1}) \cdots t'_{m+n,m+\tau(n)}(u+\ga_{m+n}).
\end{split}
\end{equation}
\ede

\bth[{\cite{Gow2}}]
The coefficients of $B(u)$ generate the center of $\Ymn$.
\eth

\bde\label{def:quasi-det}
Let $A$ be a matrix over a ring. The \emph{quasideterminant} $|A|_{ij}$ is defined (if exists) as
\[
|A|_{ij} = a_{ij} - r_i A_{ij}^{-1} c_j,
\]
where $r_i$ is the $i$-th row with $a_{ij}$ removed, $c_j$ is the $j$-th column with $a_{ij}$ removed, and $A_{ij}$ is the submatrix obtained by deleting the $i$-th row and $j$-th column.
\ede

It can also be expressed as
\beq
|A|_{ij} =
\left|
\begin{matrix}
 a_{11} & \cdots & a_{1j} & \cdots & a_{1n} \\
 \vdots & \ddots & \vdots & \ddots & \vdots \\
 a_{i1} & \cdots & \fbox{$a_{ij}$} & \cdots & a_{in}\\
 \vdots & \ddots & \vdots & \ddots & \vdots \\
 a_{n1} & \cdots & a_{nj} & \cdots & a_{nn} \\
\end{matrix}
\right|.
\eeq
If $A$ is invertible, then
\beql{quasi-det}
|A|_{ij} = \left((A^{-1})_{ji}\right)^{-1}.
\eeq

For the $(m+n)\times (m+n)$ matrix $T(u)$, let $T^{(i)}(u)$ be the submatrix consisting of its first $i$ rows and first $i$ columns.
We define $d_i(u)=|T^{(i)}(u)|_{ii}$.

\bth[{\cite{Gow1}}]\label{thm:dec} Berezinian can be expressed by quasideterminants as:
\[
B(u) = d_1(u+\ga_1) d_2(u+\ga_2) \cdots d_m(u+\ga_m) d_{m+1}(u+\ga_{m+1})^{-1} \cdots d_{m+n}(u+\ga_{m+n})^{-1}.
\]
\eth

For $1\leq i\leq m$ and $1\leq j\leq n$, let us define
\[
B_i(u) = d_1(u+\ga_1)\cdots d_i(u+\ga_i),
\]
\[
B_{m+j}(u) = d_1(u+\ga_1) d_2(u+\ga_2) \cdots d_m(u+\ga_m) d_{m+1}(u+\ga_{m+1})^{-1} \cdots d_{m+j}(u+\ga_{m+j})^{-1}.
\]
It is clear that the coefficients of $B_i(u)$ and $B_{m+j}(u)$ generate the center of $\Y(\gl_i)$ and $\Y(\gl_{m|j})$ respectively.

\bpr\label{copro on B}
For $1\leq k\leq m+n$, we have the equality
\[
\Delta^{(N-1)}(B_k(u))=B_k(u)^{\otimes N}+\text{terms of  higher degree}.
\]
\epr
\bpf
Use the root grading, in which $t_{ij}^{(p)}$ has weight $\epsilon_i-\epsilon_j$. The triangular Gauss coproduct \cite[Proposition 2.7 and Lemma 2.8]{LM} says that the diagonal term of
$\Delta(d_i(u))$ is $d_i(u)\otimes d_i(u)$; every remaining term strictly increases the first weight in the positive-root order and has the opposite weight in the second factor. Choose any total order
refining this order on the finite set of weights in question. Products, spectral shifts and inverses of series with invertible leading coefficient preserve this triangularity. The formulas of
Theorem~\ref{thm:dec} prove the assertion for $B_k(u)$. Iteration gives $N$ factors. The polynomial normalization merely multiplies each factor by known scalar series and does not change
this argument; it does not make a leading corner a Hopf subalgebra.
\epf
\subsection{Homomorphisms}
There is an injective homomorphism $\iota:U(\gl_{m|n})\rightarrow \Ymn$ given by
\[
\iota: E_{ij}\rightarrow t_{ij}^{(1)}(-1)^{\bi}.
\]
and a surjective homomorphism $\pi_{m|n}:\Ymn\rightarrow U(\gl_{m|n})$ given as follows:
\[
\pi_{m|n}: t_{ij}(u)\rightarrow \delta_{ij}+E_{ij}(-1)^{\bi}u^{-1},
\]
which is called the evaluation homomorphism.
Under the homomorphism $\iota$, we may regard $E_{ij}$ as an element of $\Ymn$.
It also gives a $\BZ^2$-grading on $\gl_{m|n}$ through this embedding.

\bpr[\cite{Naz1,Naz2}]
The map $\omega_{m|n}: \Ymn \to \Ymn$ defined by
\[
\omega_{m|n}(T(u)) = T(-u)^{-1}.
\]
is an automorphism.
\epr

\bpr
Let $f(u)$ be any formal power series in $u^{-1}$ with the leading term $1$. Then the map $\omega_f: \Ymn\to\Ymn$ defined by
\[
\omega_f(t_{ij}(u))=f(u)t_{ij}(u)
\]
is an automorphism.
\epr

Let $\phi_r : Y(\gl_{m|n}) \hookrightarrow Y(\gl_{m+r|n})$ be the natural inclusion sending $t_{ij}(u) \mapsto t_{i+r,j+r}(u)$.
Define $\psi_r : Y(\gl_{m|n}) \to Y(\gl_{m+r|n})$ by
\[
\psi_r = \omega_{m+r|n} \circ \phi_r \circ \omega_{m|n}.
\]

\bpr[\cite{Gow2}]\label{prop:psir}
For $1\leq i,j\leq m+n$, we have
\[
\psi_r (t_{ij}(u))=\left|
\begin{matrix}
 t_{11}(u) & \cdots &t_{1,r}(u) & t_{1,r+j}(u) \\
 \vdots & \ddots & \vdots & \vdots \\
 t_{r1}(u) & \cdots &t_{r,r}(u) & t_{r,r+j}(u) \\
 t_{r+i,1}(u) & \cdots &t_{r+i,r}(u) & \fbox{$t_{r+i,r+j}(u)$} \\
\end{matrix}
\right|.
\]
\epr

As a immediate consequence, there is
\bco[\cite{Gow2}]
For $1\leq k\leq m+n$, we have
\[
\psi_r(d_k(u))=d_{k+r}(u).
\]
\eco

Let us define $x_k(u), y_k(u)$ by
\[
x_k(u)=\left|
\begin{matrix}
 t_{11}(u) & \cdots &t_{1,k-1}(u) & t_{1,k+1}(u) \\
 \vdots & \ddots & \vdots & \vdots \\
 t_{k,1}(u) & \cdots &t_{k,k-1}(u) & \fbox{$t_{k,k+1}(u)$} \\
\end{matrix}
\right|,
\]
\[
y_k(u)=\left|
\begin{matrix}
 t_{11}(u) & \cdots & t_{1k}(u) \\
 \vdots & \ddots & \vdots \\
 t_{k-1,1}(u) & \cdots & t_{k-1,k}(u)\\
 t_{k+1,1}(u) & \cdots & \fbox{$t_{k+1,k}(u)$} \\
\end{matrix}
\right|.
\]
\ble\label{xy}
$[d_k(u),E_{k,k+1}]=x_k(u)$ and
$y_k(u)=(-1)^{\bar k+\overline{k+1}}[E_{k+1,k},d_k(u)]$.
\ele
\bpf
With Definition \ref{def:quasi-det} we can write
\[
d_{k}(u)=t_{kk}(u)-(t_{k1}(u),\cdots,t_{k,k-1}(u))\begin{pmatrix}&t_{11}(u)&\cdots&t_{1,k-1}(u)\\
  &\vdots &\ddots &\vdots\\
  &t_{k-1,1}(u) &\cdots &t_{k-1,k-1}(u)\\
\end{pmatrix}^{-1}\begin{pmatrix}
  &t_{1k}(u)\\
  &\vdots\\
  &t_{k-1,k}(u)\\
\end{pmatrix}
\]

For $1\leq i,j\leq k$, taking the coefficients of $v^{-1}$ in $[t_{ij}(u),t_{k,k+1}(v)]$ yields $[t_{ij}(u),E_{k,k+1}]=0$ unless $j=k$. When $j=k$,
$[t_{ik}(u),E_{k,k+1}]=t_{i,k+1}(u)$.

This implies that
\[
\begin{split}
[
&d_k(u),E_{k,k+1}]\\
=&t_{k,k+1}(u)-(t_{k1}(u),\cdots,t_{k,k-1}(u))\begin{pmatrix}&t_{11}(u)&\cdots&t_{1,k-1}(u)\\
  &\vdots &\ddots &\vdots\\
  &t_{k-1,1}(u) &\cdots &t_{k-1,k-1}(u)\\
\end{pmatrix}^{-1}\begin{pmatrix}
  &t_{1,k+1}(u)\\
  &\vdots\\
  &t_{k-1,k+1}(u)\\
\end{pmatrix}\\
=&x_k(u).
\end{split}
\]

For the second relation the RTT coefficient calculation gives
\[
[E_{k+1,k},t_{kj}(u)]=(-1)^{\bar k+\overline{k+1}}t_{k+1,j}(u)
\]
for $j\le k$, while the preceding principal block commutes with
$E_{k+1,k}$. Applying this to the same quasideterminant proves the
second formula, including its minus sign at $k=m$.
\epf

Define
\[
\phi_{k-1}(t_{ij}(u))=t_{\,i+k-1,\,j+k-1}(u).
\]
This yields a homomorphism
\[
\phi_{k-1}:
\begin{cases}
\Y(\mathfrak{gl}_2)\to \Y(\mathfrak{gl}_{k+1}), & k<m,\\
\Y(\mathfrak{gl}_{1|1})\to \Y(\mathfrak{gl}_{m|1}), & k=m,\\
\Y(\mathfrak{gl}_{0|2})\to \Y(\mathfrak{gl}_{m|\,k-m+1}), & k>m.
\end{cases}
\]
So we can define $\psi_{k-1}$ on all rank $2$ super Yangian according to $\phi_{k-1}$ defined above.
%


The proof of the following lemma is same as that of Lemma~1.11.2 in~\cite{Mol1}.


\ble\label{lem:psik}
\beq
\psi_{k-1}(t_{11}(u))=d_k(u),\quad\psi_{k-1}(t_{12}(u))=x_k(u),\quad \psi_{k-1}(t_{21}(u))=y_k(u).
\eeq
\ele


We recall another homomorphism between super-Yangians given in \cite{Gow2},
\ble[Gow \cite{Gow2}]\label{lem:rho}
The map $\rho_{m|n}:\Ymn\rightarrow \Y(\gl_{n|m})$ defined by
\[
\rho_{m|n}(t_{ij}(u))=t_{m+n+1-i,m+n+1-j}(-u)
\]
is an associative algebra isomorphism.
\ele
\ble[Nazarov \cite{Naz1,Naz2}]\label{lem:st}
The assignment $\text{st}:t_{ij}(u)\rightarrow (-1)^{\bj(\bi+1)}t_{ji}(-u)$ is an automorphism of $\Ymn$.
\ele


\subsection{Relations between Drinfeld generators}
The following lemma is useful in what follows.
\ble\label{lem:com1}
Let $1\leq k<l\leq m+n-1$, then
\begin{enumerate}
\item $(u-v)[d_k(u),y_{k}(v)]=(-1)^{\bk}(y_k(u)d_k(v)-y_k(v)d_k(u))$,
\item $(u-v)[d_k(u),x_{k}(v)]=(-1)^{\bk}(x_k(v)d_k(u)-x_k(u)d_k(v))$,
\item $y_k(u)d_k(v)-y_k(v)d_k(u)=d_k(v)y_k(u)-d_k(u)y_k(v)$,
\item $x_k(v)d_k(u)-x_k(u)d_k(v)=d_k(u)x_k(v)-d_k(v)x_k(u)$,
\item $[d_k(u),x_l(v)]=[d_k(u),y_l(v)]=0$.
\end{enumerate}
\ele
\bpf
According to Lemma \ref{lem:psik}, we can write $d_k(u)=\psi_{k-1}(t_{11}(u))$ and $y_{k}(v)=\psi_{k-1}(t_{21}(v))$,
where  $t_{11}(u),t_{21}(u)$ are elements in a space determined by the parameter $k$.
Then according to \eqref{def-rel}, we obtain Equation $(1)$.
Equation
$(2)$ can be proved  by the same argument.

The identity
\[
t_{21}(u)t_{11}(v)-t_{21}(v)t_{11}(u)=t_{11}(v)t_{21}(u)-t_{11}(u)t_{21}(v),
\]
 holds in $\Y(\gl_{a|b})$ for any nonnegative $a,b$. Applying $\psi_{k-1}$ to  this identity we get Equation $(3)$.
 Following the same method, we can derive
Equation $(4)$.

According to Lemma \ref{xy},
we write $x_l(v)$ as $[d_l(v),E_{l,l+1}]$
Then,
\[
[d_k(u),x_l(v)]=[[d_k(u),d_l(v)],E_{l,l+1}]+[d_l(v),[d_k(u),E_{l,l+1}]].
\]
As the method shown in the proof of Lemma \ref{xy}, we can see that $[d_k(u),E_{l,l+1}]=0$ when $k<l$.
As a result, there is $[d_k(u),x_l(v)]=0$. Similarly, $[d_k(u),y_l(v)]=0$.

This proves $(5)$.

\epf

\section{The skew representation}

\subsection{Gelfand-Tsetlin basis}\label{sec:GT bas}

We have already known that given a covariant weight $\la$, there exists a unique simple finite dimensional $\gl_{m|n}$-module $L(\la)$.
It is equipped with a basis parameterized by the Gelfand-Tsetlin patterns $\La=(\la_{ki}\in\BZ|1\leq i\leq k\leq m+n)$,
which satisfy the following conditions:
\begin{enumerate}
 \item $\lambda_{m+n,i}=\lambda_i, \quad 1\leq i\leq m+n$;
 \item  $\la_{k,i}-\la_{k-1,i}\equiv\theta_{k-1,i}\in\{0,1\}, 1\leq i\leq m;m+1\leq k\leq m+n$;
 \item $\la_{ki}-\la_{k,i+1}\in \mathrm{Z}_{\geq 0}, 1\leq i\leq m-1;m+1\leq k\leq m+n-1$;
 \item  $\la_{k+1,i}-\la_{ki}\in \mathrm{Z}_{\geq 0}$ {\it and} $\la_{k,i}-\la_{k+1,i+1}\in \mathrm{Z}_{\geq 0},$
$1\leq i\leq k\leq m-1$ \text{ or } $m+1\leq i\leq k\leq m+n-1$;
  \item
$m+1\le k\le m+n$:
$\lambda_{km}\ge \#\{i:\lambda_{ki}>0,\; m+1\le i\le k\};$
  \item
if $\lambda_{ m+1,m}=0$, then $\theta_{mm}=0$.
\end{enumerate}
We denote the set of all patterns satisfying the above condition by $\mathscr{S}_{\la}$. It can be written in the triangular form
\begin{equation*}
\resizebox{.9\textwidth}{!}{$\begin{array}{cccccccc}
   \la_{m+n,1}   &  \cdots      & \la_{m+n,m} & \la_{m+n,m+1} & \cdots & \la_{m+n,m+n-1} & \la_{m+n,m+n} \\
  \la_{m+n-1,1} &  \cdots     & \la_{m+n-1,m} & \la_{m+n-1,m+1} & \cdots & \la_{m+n-1,m+n-1} & \\
  \vdots    &  \vdots & \vdots & \vdots   & \reflectbox{$\ddots$}    \\
  \la_{m+1,1} &  \cdots  & \la_{m+1,m} & \la_{m+1,m+1}    \\
  \la_{m,1} &  \cdots  & \la_{m,m}                           \\
  \la_{m-1,1} &  \cdots                             \\
  \vdots  &  \reflectbox{$\ddots$}                              \\
  \la_{11}\\
\end{array}$}
\end{equation*}

A certain class of finte dimensional \emph{essentially typical} representations of the Lie superalgebra $\gl(m|n)$ were constructed in \cite{Palev1989a},
\cite{Palev1989b}.
A Gelfand-Tsetlin basis for these representations was obtained therein.
This construction was later generalized to \emph{covariant tensor} modules in \cite{SV}.
Another explicit construction of covariant tensor modules, based on super Young tableaux, was given by Molev \cite{Mol2}.
In the following we use the formulas from \cite{FSZ},
which are obtained by modifying the formulas of Palev \cite{Palev1989b}, Stoilova and Van der Jeugt \cite{SV}.

\bth[\cite{Mol2} \cite{SV}]\label{for thm}
There exists a basis $\{\xi_{\Lambda}\}$ in $L(\lambda)$ parametrized by all Gelfand-Tsetlin patterns $\Lambda$, the action of generators of $\gl_{m|n}$ is given by the formulas
\beq
E_{kk}\xi_{\La}=\left(\sum_{i=1}^{k}\la_{ki}-\sum_{j=1}^{k-1}\la_{k-1,j}\right)\xi_{\La},\quad 1\leq k\leq m+n;
\eeq
\beq
E_{k,k+1}\xi_{\La}=-\sum_{i=1}^{k}\frac{\Pi_{j=1}^{k+1}(l_{k+1,j}-l_{ki}) }
  {\Pi_{j\neq i,j=1}^{k} (l_{kj}-l_{ki}) }\xi_{\La+\delta_{ki}},
\quad1\leq k\leq m-1;
\eeq
\beq
E_{k+1,k}\xi_{\La}=\sum_{i=1}^{k}\frac{\Pi_{j=1}^{k-1}(l_{k-1,j}-l_{ki})}
{\Pi_{j\neq i,j=1}^{k}(l_{kj}-l_{ki})}\xi_{\La-\delta_{ki}},
\quad  1\leq k\leq m-1;
\eeq
\beq
\begin{split}
E_{m,m+1}\xi_{\La}&=\sum_{i=1}^{m}\theta_{mi}(-1)^{i-1}(-1)^{\theta_{m1}+\ldots+\theta_{m,i-1}}\\
&\quad\times  \frac{\Pi_{1\leq j< i} (l_{mj}-l_{mi}-1)}
  {\Pi_{i<j\leq m} (l_{mj}-l_{mi})
    \Pi_{j\neq i,j=1}^{m}(l_{m+1,j}-l_{mi}-1)}
    \xi_{\La+\delta_{mi}},
\end{split}
\eeq
\beq
\begin{split}
&E_{m+1,m} \xi_{\La}=\sum_{i=1}^{m}(1-\theta_{mi})(-1)^{i-1}(-1)^{\theta_{m1}+\ldots+\theta_{m,i-1}}\\
& \times  \frac{(l_{m,i}-l_{m+1,m+1})\Pi_{ i<j\leq m} (l_{mj}-l_{mi}+1)\Pi_{j=1}^{m-1}(l_{m-1,j}-l_{mi})}
  {\Pi_{1\leq j< i} (l_{mj}-l_{mi})} \xi_{\La-\delta_{mi}},
\end{split}
\eeq
\beq
\begin{split}
E_{k,k+1} \xi_{\La}&=\sum_{i=1}^{m}\theta_{ki}(-1)^{\theta_{k1}+\ldots+\theta_{k,i-1}
+\theta_{k-1,i+1}+\ldots+\theta_{k-1,m}}(1-\theta_{k-1,i})\\
&\qquad\times
\prod_{j\neq i,j =1}^{m}\left(\frac{l_{kj}-l_{ki}-1}{l_{k+1,j}-l_{ki}-1}\right)
\xi_{\La+\delta_{ki}}\\
&\quad-\sum_{i=m+1}^{k}
\Pi_{j=1}^{m}\left(\frac{(l_{kj}-l_{ki})(l_{kj}-l_{ki}+1)}{(l_{k+1,j}-l_{ki})(l_{k-1,j}-l_{ki}+1)} \right)\\
&\qquad\times\frac{\Pi_{j=m+1}^{k+1}(l_{k+1,j}-l_{ki})}
{\Pi_{j\neq i,j=m+1}^{k} (l_{kj}-l_{ki})}\xi_{\La+\delta_{ki}} ,
\quad  m+1\leq k\leq m+n-1;
\end{split}
\eeq
\beq
\begin{split}
E_{k+1,k}\xi_{\La}&=
\sum_{i=1}^{m}\theta_{k-1,i}(-1)^{\theta_{k1}+\ldots+\theta_{k,i-1}+\theta_{k-1,i+1}+\ldots+\theta_{k-1,m}}(1-\theta_{ki})
\\
&\qquad \times
\prod_{j\neq i,j=1}^{m}\left(\frac{l_{kj}-l_{ki}+1}{l_{k-1,j}-l_{ki}+1}\right)\\
&\qquad\times\frac{\Pi_{j=m+1}^{k+1}(l_{k+1,j}-l_{ki})\Pi_{j=m+1}^{k-1}(l_{k-1,j}-l_{ki}+1)}
  {\Pi_{j=m+1}^{k} (l_{kj}-l_{ki})(l_{kj}-l_{ki}+1)}\xi_{\La-\delta_{ki}}
  \\
&\quad  + \sum_{i=m+1}^{k}
\frac{\prod_{j=m+1}^{k-1}(l_{k-1,j}-l_{ki})}{\prod_{j\neq i,j =m+1}^{k}(l_{k,j}-l_{ki})}
\xi_{\La-\delta_{ki}} \quad
 m+1\leq k\leq m+n-1.
\end{split}
\eeq
where
\begin{equation*}
l_{ki}=\la_{ki}-i+1, (1\leq i\leq m); \quad l_{kj}=-\la_{kj}+j-2m, (m+1\leq j\leq k).
\end{equation*}
The arrays $\La\pm \delta_{ki}$ are obtained from $\La$ by replacing $\la_{ki}$ by $\la_{ki}\pm1$. We assume that
$\xi_{\La}=0$ if the array $\La$ is not a Gelfand-Tsetlin pattern.
\eth

It is clear that every $\xiL$ is an eigenvector of $t_{ii}(u)$ with eigenvalue,
\[
\la_i(u)=u+(-1)^{\bi}(\sum_{j=1}^i\la_{ij}-\sum_{j=1}^{i-1}\la_{i-1,j}).
\]

For a Gelfand-Tsetlin pattern $\La$,
we give $\xiL$ a $\BZ^2$-grading by
\[
\deg (\xiL)=(\sum_{k=1}^m\sum_{i=1}^k \la_{ki},\sum_{k=m+1}^{m+n}\sum_{i=1}^k \la_{ki}).
\]
It can be seen that $\deg(x.v)=\deg(x)+\deg(v)$ for a homogeneous element $x\in \gl_{m|n}$ and a homogeneous element $v\in V(\la)$. According to the
  partial order on $\BZ^2$ given at the end of Section \ref{sec:def}, we also define a partial order on homogeneous elements in $L(\la)$.
Define $v< w$ if $\deg(v)<\deg(w)$.

\subsection{Skew representations}\label{sec:skew mod}
In the remaining of the paper, we may use $ut_{ij}(u)$ to instead $t_{ij}(u)$ to avoid some factors in the denominators. And we still denote it by $t_{ij}(u)$.

Now consider the simple finite dimensional $\gl_{m+r|n}$-module $L(\la)$ with $\la=(\la_1,\ldots,\la_{m+r+n})$.
Let $\mu=(\mu_1,\ldots,\mu_r)$ satisfy that
\beql{admi}
\la_{m+i}\leq \mu_i\leq \la_i.
\eeq
Define $L(\la\slash\mu)$ as the subspace of $L(\lambda)$ formed by all singular vectors with respect to $\gl_{r}$ of weight $\mu$, i.e.,
\[
L(\la\slash\mu)=\{v\in V(\lambda)|E_{ii}v=\mu_iv,E_{jk}v=0\text{ for }1\leq i\leq r,1\leq j<k\leq r\}.
\]
Then $L(\la\slash\mu)$ takes a basis parameterized by the patterns
$\La=(\la_{ij}|1\leq j\leq i\leq m+r+n)$ such that $\la_{kj}=\mu_j$ for $k=1,\ldots,r$.
Denote the set consisting of all such patterns by $\mathscr{S}_{\la,\mu}$.

It is known that $\pi_{m+r|n}\psi_r$ commute with the subalgebra $\gl_r$,
which implies that $L(\la\slash\mu)$ can be regarded as a $\Ymn$-module (called the skew module) through the homomorphism $\pi_{m+r|n}\psi_r$.

The relation \eqref{RTT} implies that for any $h\in \BC$ the assignment $t_{ij}(u)\rightarrow t_{ij}(u+h)$ determines an automorphism of the algebra $\Ymn$.
We will denote by $L_{h}(\la\slash\mu)$ the $\Ymn$-module obtained from $L(\la\slash\mu)$ by the pullback through this automorphism.

Given $\La\in\mathscr{S}_{\la,\mu}$, define
\[
\chil(u)=\frac{\prod_{k=1}^{r+i}(u+r+l_{r+i,k})}{\prod_{k=1}^{r}(u+r+o_k)},
\]
and
\[
\chjl(u)=\frac{\prod_{k=1}^{r+m}(u+r+l_{r+m+j,k})}{\prod_{k=1}^{r}(u+r+o_k)\prod_{k=1}^{j}(u+r+l_{r+m+j,r+m+k})},
\]
for $1\leq i\leq m$ and $1\leq j\leq n$, where $o_k=\mu_k-k+1$ for $1\leq k\leq r$.
\ble\label{lem:chkL}
Let $\La$ be a Gelfand-Tsetlin pattern in $\mathscr{S}_{\la,\mu}$, i.e., $\xiL\in L(\la\slash\mu)$.
Then we have
\begin{enumerate}
\item $B_i(u)\xiL=\chil(u)\xiL$,
\item $B_{m+j}(u)\xiL=\chjl(u)\xiL$.
\end{enumerate}
\ele

Let
\[
\zeil(u+\ga_i)=\frac{\prod_{k=1}^{r+i}(u+l_{r+i,k}+r)}{\prod_{k=1}^{r+i-1}(u+l_{r+i-1,k}+r)},
\]
\[
\zejl(u+\ga_{m+j})=\frac{\prod_{k=1}^{m+r}(u+l_{m+r+j-1,k}+r)}{\prod_{k=1}^{m+r}(u+l_{m+r+j,k}+r)}\frac{\prod_{k'=1}^{j}(u+r+l_{m+r+j,m+r+k'})}{\prod_{k'=1}^{j-1}(u+r+l_{m+r+j-1,m+r+k'})}.
\]
\ble\label{lem:zek}
Let $\La$ be a Gelfand-Tsetlin pattern in $\mathscr{S}_{\la,\mu}$, we have
\begin{enumerate}
\item $d_i(u)\xiL=\zeil(u) \xiL$.
\item $d_{m+j}(u)\xiL=\zejl(u) \xiL$.
\end{enumerate}
\ele

For every $s=1,\ldots,N$ fix some $\hs\in\BC$ along with a pair of sequences of intergers
\[
\la^{(s)}=(\la_1^{(s)},\ldots,\la_{r^{(s)}}^{(s)},\la_{r^{(s)}+1}^{(s)},\ldots,\la_{m+r^{(s)}}^{(s)},\la_{m+r^{(s)}+1}^{(s)},\ldots,\la_{m+r^{(s)}+n}^{(s)})
\]
and
\[
\mu^{(s)}=(\mu_1^{(s)},\ldots,\mu_{r^{(s)}}^{(s)}).
\]
Here $r^{(s)}$ is a nonnegative integer, $\la^{(s)}$ is a covariant weight and $\mu^{(s)}$ satisfies \eqref{admi} for every $s$.

Through the comultiplication \eqref{copro}, we regard the following as a $\Ymn$-module,
\[
W=L_{h^{(1)}}(\la^{(1)}\slash\mu^{(1)})\otimes\cdots\otimes L_{h^{(N)}}(\la^{(N)}\slash\mu^{(N)}).
\]

\section{Gelfand-Tsetlin basis for tensor product of skew modules}\label{sec:irr}
In this section, we give a Gelfand-Tsetlin type basis of
\[
W=L_{h^{(1)}}(\la^{(1)}\slash\mu^{(1)})\otimes\cdots\otimes L_{h^{(N)}}(\la^{(N)}\slash\mu^{(N)}).
\]
when $h^{(s)}- h^{(t)}\notin \BZ$ for all $s\neq t$ which is analogous to \cite{Mol,NT1,NT2}.

Put
\beq\label{eq:skew_chi_product}
\chkll(u)=\prod_{s=1}^N\chi_{k,\La^{(s)}}(u+h^{(s)}),\qquad
\zekll(u)=\prod_{s=1}^N\zeta_{k,\La^{(s)}}(u+h^{(s)}).
\eeq
We use the polynomial normalization allowed in
Section~2.3.

For any $\La^{(s)},\Omega^{(s)}\in
\mathscr S_{\la^{(s)},\mu^{(s)}}$, we have
\beq\label{eq:skew_separation}
\chi_{k,\La^{(1)},\ldots,\La^{(N)}}(u)
=\chi_{k,\Omega^{(1)},\ldots,\Omega^{(N)}}(u)
\ \ (1\le k\le m+n)
\quad\Longrightarrow\quad
\La^{(s)}=\Omega^{(s)}\ \ (1\le s\le N).
\eeq
Indeed, roots and poles belonging to distinct integral cosets cannot cancel. Within a fixed factor, the formula in Lemma~\ref{lem:chkL} and the Gelfand--Tsetlin inequalities recover the successive rows.

Let $\La=(\La^{(1)},\ldots,\La^{(N)})$ be a collection of pattern and
\[
v_{\La}=\xi_{\La^{(1)}}\otimes\cdots\otimes\xi_{\La^{(N)}}.
\]
There are finitely many coefficients of the series $B_k(u)=\sum_{p\in\BZ}B_k^{(p)}u^{-p}$ and complex numbers $c_{k,p}$ such that
\beq\label{eq:skew_A}
A=\sum_{k,p}c_{k,p}B_k^{(p)},
\eeq
has pairwise distinct diagonal entries $a_\La$ on the vectors $v_\La$. Then set
\beql{eigvec}
\xi_\La=P_\La v_\La,
\eeq
where $P_\La=\prod_{\Omega\neq\La}\frac{A-a_\Omega}{a_\La-a_\Omega}$.

\bpr\label{thm:bas}
Let $\La=(\La^{(1)},\ldots,\La^{(N)})$, where $\La^{(s)}\in\mathscr S_{\la^{(s)},\mu^{(s)}}\ (1\le s\le N)$ and $\xi_\La$ be defiend as above, then
\[
\{\xi_{\La}=\xiLL:\
\La^{(s)}\in\mathscr S_{\la^{(s)},\mu^{(s)}},\ 1\le s\le N\}
\]
is a basis of $W$. Moreover, for $1\le k\le m+n$,
\beq\label{eq:skew_eigen_formula}
d_k(u+\ga_k)\xi_{\La}=\zekll(u+\ga_k)\xi_{\La}.
\eeq
\epr
\bpf
By Proposition~\ref{copro on B} and Lemma~\ref{lem:chkL},
\beq\label{eq:skew_triangular}
B_k(u)v_\La
=\chi_{k,\La}(u)v_\La
+\sum_{\Omega>\La}b_{k,\La,\Omega}(u)v_\Omega,
\qquad 1\le k\le m+n.
\eeq
Equation~\eqref{eq:skew_separation}
implies that finitely many coefficients distinguish the diagonal
characters. Hence the $c_{k,p}$ in~\eqref{eq:skew_A} may be
chosen outside finitely many hyperplanes so that
$a_\La\ne a_\Omega$ whenever $\La\ne\Omega$. Consequently,
\beq\label{eq:skew_leading}
P_\La v_\La
=v_\La+\sum_{\Omega>\La}c_{\La,\Omega}v_\Omega\ne0.
\eeq
Since $P_\La P_\Omega=\delta_{\La\Omega}P_\La$ on $W$, the vectors in~\eqref{eigvec} form a basis. By
$[d_k(u),B_l(v)]=0$, the one-dimensional image $P_\La W$ is $d_k(u)$-stable.
Comparing its diagonal entry in \eqref{eq:skew_triangular} and applying Theorem~\ref{thm:dec}
gives~\eqref{eq:skew_eigen_formula}.
\epf

Let $\La_0^{(s)}$ be the possible highest pattern with the $r^{(s)}$-row fixed as shown in Section~\ref{sec:GT bas} and put $\xi_0=\xi_{\La_0^{(1)}}\otimes\cdots\otimes\xi_{\La_0^{(N)}}$.
\ble\label{lem:stp}
For $1\le k\le m+n$,
\[
d_k(u+\ga_k)\xi_0=\zeta_{k,\La_0^{(1)},\ldots,\La_0^{(N)}}(u+\ga_k)\xi_0.
\]
\ele
\bpf
Note that for the possible highest pattern $\La_0$ with the $r$-row fixed, $\xi_{\La_0}$ lies in finite dimensional simple $\gl_{m+r,n}$-module $V(\la)$,
according to Definition \ref{def:quasi-det} and Proposition \ref{prop:psir}, there is
\[
\begin{split}
&\psi_{r}(t_{ij}(u))\xi_{\La_0}\\
=&t_{r+i,r+j}(u)\xi_{\La_0}-(t_{r+i,1}(u),\cdots,t_{r+i,r}(u))\begin{pmatrix}&t_{11}(u)&\cdots&t_{1r}(u)\\
&\vdots &\ddots &\vdots\\
 &t_{r,1}(u) &\cdots &t_{rr}(u)\\
\end{pmatrix}^{-1}
\begin{pmatrix}
&t_{1,r+j}(u)\\
&\vdots\\
&t_{r,r+j}(u)\\
\end{pmatrix}\xi_{\La_0}\\
=&t_{r+i,r+j}(u)\xi_{\La_0}.
\end{split}
\]
It implies that
\[
\begin{split}
t_{kk}(u)\xi_0=&t_{k+r^{(1)},k+r^{(1)}}(u+h^{(1)})\xi_{\La_0^{(1)}}\otimes\cdots\otimes t_{k+r^{(N)},k+r^{(N)}}(u+h^{(N)})\xi_{\La_0^{(N)}}\\
=&(u+h^{(1)}+(-1)^{\overline{k+r^{(1)}}}\la_{k+r^{(1)}})\cdots (u+h^{(N)}+(-1)^{\overline{k+r^{(N)}}}\la_{k+r^{(N)}})\xi_0\\
=&\ze_{k,\La_0^{(1)},\ldots,\La_0^{(N)}}(u)\xi_0.
\end{split}
\]
It also implies that $t_{ik}(u)\xi_0=0$ if $i\leq k$.

Let us consider the equality $(T^{(k)}(u))^{-1}T^{(k)}(u)=1$. It implies that
\[
\begin{split}
\xi_0=&\sum_{i=1}^{k}\left((T^{(k)}(u))^{-1}\right)_{ki}t_{ik}(u)\xi_0\\
=&\left((T^{(k)}(u))^{-1}\right)_{kk}t_{kk}(u)\xi_0.
\end{split}
\]
But we have already seen $\left((T^{(k)}(u))^{-1}\right)_{kk}=d_k(u)^{-1}$ due to \eqref{quasi-det}. Thus, we finish our proof.
\epf

We calculate the raising and lowering operators on the basis \eqref{eigvec}.  Put
\[
e_k(u)=d_k(u)^{-1}x_k(u),\qquad f_k(u)=y_k(u)d_k(u)^{-1},\qquad \phi_k(u)=d_k(u)^{-1}d_{k+1}(u).
\]
Write $e_k(u)=\sum_{p\ge1}e_k^{(p)}u^{-p}$ and similarly for $f_k(u)$.  Since $d_k(u)=u^N+O(u^{N-1})$ on $W$, the leading coefficients of $x_k(u)$ and $y_k(u)$ are coefficients of $u$ of degree $N-1$; in the present section we denote those coefficients by $x_k^{(1)}=[u^{N-1}]x_k(u)$ and
$y_k^{(1)}=[u^{N-1}]y_k(u)$.
By the evaluation homomorphism and Lemma~\ref{xy},
\beq\label{eq:skew_first_coefficients}
\begin{aligned}
e_k^{(1)}=x_k^{(1)}&=(-1)^{\bk}E_{k,k+1},\\
f_k^{(1)}=y_k^{(1)}&=(-1)^{\overline{k+1}}E_{k+1,k}.
\end{aligned}
\eeq
The coefficients $E_{k,k+1}$ and $E_{k+1,k}$ on each tensor factor are the fractions in Theorem~\ref{for thm}, with the row index $k$ replaced by $k+r^{(s)}$.

The Gauss decomposition and the rank-two RTT relations, transported by $\psi_{k-1}$, give
\beq\label{eq:skew_ef_comm}
\begin{aligned}
[e_k^{(p)},f_l^{(q)}]&=-\delta_{kl}(-1)^{\overline{k+1}}\phi_k^{(p+q-1)},\\
[e_k(u),f_k^{(1)}]&=(-1)^{\overline{k+1}}(1-\phi_k(u)).
\end{aligned}
\eeq
The brackets are supercommutators.  In particular, $[B_j(u),e_k(v)]=[B_j(u),f_k(v)]=0$ for $j\ne k$.

The relations of Lemma~\ref{lem:com1}, together with the $B_j$-characters, imply that the coefficient of $\xi_\Omega$ in $e_k(u)\xi_\La$ (respectively $f_k(u)\xi_\La$) vanishes unless $\Omega^{(s)}=\La^{(s)}+\delta_{k+r^{(s)},i}\in\mathscr S_{\la^{(s)},\mu^{(s)}}$ (respectively $\Omega^{(s)}=\La^{(s)}-\delta_{k+r^{(s)},i}\in\mathscr S_{\la^{(s)},\mu^{(s)}}$)
for one $s,i$, with all other components unchanged.
Indeed, suppose that the coefficient of \(\xi_\Omega\) in \(x_k(u)\xi_\Lambda\) is nonzero. By the commutation relations, we have
\[
\chi_{j,\Omega}(u)=\chi_{j,\Lambda}(u),\qquad j\ne k.
\]
Moreover, the matrix-coefficient formula in \cite[Proposition 3.4]{Lu} shows that the corresponding coefficient of \(d_k(u)^{-1}x_k(u)\) has the form \(C/(u-a)\) for some \(C\ne0\). The pole belongs to exactly one coset
$-h^{(s)}+\BZ$. Separating the other cosets in the character ratio leaves their factors unchanged by \eqref{eq:skew_separation}.
In factor $s$, the single-factor character calculation identifies the change as one entry in row $k+r^{(s)}$, with the indicated sign.

Let $\Omega$ differ from $\La$ only by $\Omega^{(s)}=\La^{(s)}+\delta_{k+r^{(s)},i}
\in\mathscr S_{\la^{(s)},\mu^{(s)}}$.
In the already fixed normalization~\eqref{eigvec}, set
\beq\label{eq:skew_raise_coefficient}
C^+_{k;\Omega,\La}=[v_\Omega]P_\Omega x_k^{(1)}P_\La v_\La,
\eeq
the coefficient of the basis vector $v_\Omega$ in $P_\Omega x_k^{(1)}P_\La v_\La$ and, if $\Omega$ differs from $\La$ only by $\Omega^{(s)}=\La^{(s)}-\delta_{k+r^{(s)},i}\in\mathscr S_{\la^{(s)},\mu^{(s)}}$, put
\beq\label{eq:skew_lower_coefficient}
C^-_{k;\Omega,\La}=[v_\Omega]P_\Omega y_k^{(1)}P_\La v_\La.
\eeq
Both numbers are finite expressions in the coefficients of
Theorem~\ref{for thm} and the products~\eqref{eigvec}.

For each $s$, the ratio of the $d_k(u)$-eigenvalue under an admissable entry change is given by the following rational function, where $z=u+h^{(s)}+r^{(s)}-\gamma_k$:
\beq\label{eq:skew_ratio_table}
\begin{aligned}
&\frac{\zeta_{k,\La^{(s)}+\delta_{k+r^{(s)},i}}(u+h^{(s)})}
{\zeta_{k,\La^{(s)}}(u+h^{(s)})}
=\begin{cases}
\dfrac{z+l_{k+r^{(s)},i}^{(s)}+1}{z+l_{k+r^{(s)},i}^{(s)}},&k\le m,\\[4pt]
\dfrac{z+l_{k+r^{(s)},i}^{(s)}}{z+l_{k+r^{(s)},i}^{(s)}+1},&k>m,\ i\le m+r^{(s)},\\[4pt]
\dfrac{z+l_{k+r^{(s)},i}^{(s)}-1}{z+l_{k+r^{(s)},i}^{(s)}},&k>m,\ i>m+r^{(s)},
\end{cases}\\[5pt]
&\frac{\zeta_{k,\La^{(s)}-\delta_{k+r^{(s)},i}}(u+h^{(s)})}
{\zeta_{k,\La^{(s)}}(u+h^{(s)})}
=\begin{cases}
\dfrac{z+l_{k+r^{(s)},i}^{(s)}-1}{z+l_{k+r^{(s)},i}^{(s)}},&k\le m,\\[4pt]
\dfrac{z+l_{k+r^{(s)},i}^{(s)}}{z+l_{k+r^{(s)},i}^{(s)}-1},&k>m,\ i\le m+r^{(s)},\\[4pt]
\dfrac{z+l_{k+r^{(s)},i}^{(s)}+1}{z+l_{k+r^{(s)},i}^{(s)}},&k>m,\ i>m+r^{(s)}.
\end{cases}
\end{aligned}
\eeq
Here the first (respectively second) fraction is considered
only if $\La^{(s)}+\delta_{k+r^{(s)},i}\in
\mathscr S_{\la^{(s)},\mu^{(s)}}$ (respectively
$\La^{(s)}-\delta_{k+r^{(s)},i}\in\mathscr S_{\la^{(s)},\mu^{(s)}}$).
Taking matrix coefficients in Lemma~\ref{lem:com1}(1),(2)
gives
\[
\begin{aligned}
e_k(v)\xi_\La
&=\sum_{\substack{1\le s\le N,\ 1\le i\le k+r^{(s)}\\
\La^{(s)}+\delta_{k+r^{(s)},i}\in\mathscr S_{\la^{(s)},\mu^{(s)}}}}
\frac{C^+_{k;\La+\delta_{k+r^{(s)},i},\La}}
{v-\ga_k+h^{(s)}+r^{(s)}+c^+_{k,s,i}}
\xi_{\La+\delta_{k+r^{(s)},i}},\\[4pt]
f_k(v)\xi_\La
&=\sum_{\substack{1\le s\le N,\ 1\le i\le k+r^{(s)}\\
\La^{(s)}-\delta_{k+r^{(s)},i}\in\mathscr S_{\la^{(s)},\mu^{(s)}}}}
\frac{C^-_{k;\La-\delta_{k+r^{(s)},i},\La}}
{v-\ga_k+h^{(s)}+r^{(s)}+c^-_{k,s,i}}
\xi_{\La-\delta_{k+r^{(s)},i}}.
\end{aligned}
\]
As a result, we have
\beq\label{eq:skew_x_y_action}
\begin{aligned}
x_k(v)\xi_\La
&=\sum_{\Omega=\La+\delta_{k+r^{(s)},i}}
\frac{C^+_{k;\Omega,\La}\,\zeta_{k,\Omega}(v)}
{v-\ga_k+h^{(s)}+r^{(s)}+c^+_{k,s,i}}\xi_\Omega,\\[3pt]
y_k(v)\xi_\La
&=\zeta_{k,\La}(v)
\sum_{\Omega=\La-\delta_{k+r^{(s)},i}}
\frac{C^-_{k;\Omega,\La}}
{v-\ga_k+h^{(s)}+r^{(s)}+c^-_{k,s,i}}\xi_\Omega.
\end{aligned}
\eeq
$\Omega$ differs from $\La$ only in
one factor $s$, and $\Omega^{(s)}\in
\mathscr S_{\la^{(s)},\mu^{(s)}}$. The integer shifts are
\beq\label{eq:skew_pole_cases}
\begin{array}{c|cc}
&c^+_{k,s,i}&c^-_{k,s,i}\\ \hline
k\le m&l_{k+r^{(s)},i}^{(s)}+1&l_{k+r^{(s)},i}^{(s)}\\
k>m,\ i\le m+r^{(s)}&l_{k+r^{(s)},i}^{(s)}&l_{k+r^{(s)},i}^{(s)}-1\\
k>m,\ i>m+r^{(s)}&l_{k+r^{(s)},i}^{(s)}-1&l_{k+r^{(s)},i}^{(s)}.
\end{array}
\eeq

\bth\label{thm:ksi}
The formulas~\eqref{eq:skew_eigen_formula} and \eqref{eq:skew_x_y_action} give the actions of
$d_k(u)$, $x_k(u)$ and $y_k(u)$ on the vectors~\eqref{eigvec}.
Moreover, every coefficient $C^+_{k;\Omega,\La}$ and $C^-_{k;\Omega,\La}$ corresponding to an entry change with $\Omega^{(s)}\in\mathscr S_{\la^{(s)},\mu^{(s)}}$ is nonzero.
\eth
\bpf
It remains to prove the last assertion. By the formulas of Lemma~\ref{lem:zek}, on the vector $\xi_\La$ the series $\phi_k(v)$ acts by
\beq\label{eq:skew_phi_product}
\Phi_{k,\La}(v)=\frac{\zeta_{k+1,\La}(v)}{\zeta_{k,\La}(v)}=\prod_{s=1}^{N}R_{k,\La^{(s)}}(v+h^{(s)}+r^{(s)}-\ga_k),
\eeq
where the rational function $R_{k,\La^{(s)}}(z)$ is given explicitly as follows.
\beq\label{eq:skew_phi_even}
R_{k,\La^{(s)}}(z)=\frac{\displaystyle\prod_{j=1}^{k+r^{(s)}+1}(z+l^{(s)}_{k+r^{(s)}+1,j}+1)
\prod_{j=1}^{k+r^{(s)}-1}(z+l^{(s)}_{k+r^{(s)}-1,j})}{\displaystyle\prod_{j=1}^{k+r^{(s)}}(z+l^{(s)}_{k+r^{(s)},j})(z+l^{(s)}_{k+r^{(s)},j}+1)},\quad k<m,
\eeq
\beq\label{eq:skew_phi_mixed}
R_{m,\La^{(s)}}(z)=\frac{(z+l^{(s)}_{m+r^{(s)}+1,m+r^{(s)}+1})\displaystyle\prod_{j=1}^{m+r^{(s)}-1}(z+l^{(s)}_{m+r^{(s)}-1,j})}{\displaystyle\prod_{j=1}^{m+r^{(s)}}(z+l^{(s)}_{m+r^{(s)}+1,j})},
\eeq
\beq\label{eq:skew_phi_odd}
\begin{aligned}
R_{k,\La^{(s)}}(z)={}&\prod_{j=1}^{m+r^{(s)}}\frac{(z+l^{(s)}_{k+r^{(s)},j}-1)(z+l^{(s)}_{k+r^{(s)},j})}
{(z+l^{(s)}_{k+r^{(s)}+1,j}-1)(z+l^{(s)}_{k+r^{(s)}-1,j})}\\
&\times\frac{\displaystyle\prod_{j=m+r^{(s)}+1}^{k+r^{(s)}+1}(z+l^{(s)}_{k+r^{(s)}+1,j}-1)
\prod_{j=m+r^{(s)}+1}^{k+r^{(s)}-1}(z+l^{(s)}_{k+r^{(s)}-1,j})}{\displaystyle\prod_{j=m+r^{(s)}+1}^{k+r^{(s)}}(z+l^{(s)}_{k+r^{(s)},j}-1)(z+l^{(s)}_{k+r^{(s)},j})},\qquad k>m.
\end{aligned}
\eeq

Let $\La^{(s)}\pm\delta_{k+r^{(s)},i}\in\mathscr S_{\lambda^{(s)},\mu^{(s)}}$, and put $c=c^\pm_{k,s,i}$ as in \eqref{eq:skew_pole_cases}. To determine the order of $R_{k,\Lambda^{(s)}}(z)$ at $z=-c$, we use the following notation: $\delta_a=1$ if $a$ is true, and $0$ otherwise.

The following is a direct zero--pole count for $R_{k,\La^{(s)}}(z)$.
For $k<m$ and $k=m$, the orders at $z=-c$ are respectively
\beq\label{eq:skew_order_evenroot}
\begin{aligned}
&\sum_{j=1}^{k+r^{(s)}+1}\de_{l_{k+r^{(s)}+1,j}+1=c}+\sum_{j=1}^{k+r^{(s)}-1}
\de_{l_{k+r^{(s)}-1,j}=c}\\
&\quad-\sum_{j=1}^{k+r^{(s)}}\left(\de_{l_{k+r^{(s)},j}=c}+\de_{l_{k+r^{(s)},j}+1=c}
\right)=-1,\\
&\de_{l_{m+r^{(s)}+1,m+r^{(s)}+1}=c}+\sum_{j=1}^{m+r^{(s)}-1}\de_{l_{m+r^{(s)}-1,j}=c}\\
&\quad-\sum_{j=1}^{m+r^{(s)}}\de_{l_{m+r^{(s)}+1,j}=c}=-1.
\end{aligned}
\eeq
For $k>m$, the orders of the two factors on the right of \eqref{eq:skew_phi_odd} at $z=-c$ are respectively
\beq\label{eq:skew_order_even}
\sum_{j=1}^{m+r^{(s)}}\bigl(\de_{l_{k+r^{(s)},j}-1=c}+\de_{l_{k+r^{(s)},j}=c}
-\de_{l_{k+r^{(s)}+1,j}-1=c}-\de_{l_{k+r^{(s)}-1,j}=c}
\bigr)
\eeq
and
\beq\label{eq:skew_order_odd}
\begin{aligned}
&\sum_{j=m+r^{(s)}+1}^{k+r^{(s)}+1}\de_{l_{k+r^{(s)}+1,j}-1=c}+\sum_{j=m+r^{(s)}+1}^{k+r^{(s)}-1}
\de_{l_{k+r^{(s)}-1,j}=c}\\
&\hspace{9mm}
-\sum_{j=m+r^{(s)}+1}^{k+r^{(s)}}\bigl(\de_{l_{k+r^{(s)},j}-1=c}+\de_{l_{k+r^{(s)},j}=c}
\bigr).
\end{aligned}
\eeq

The exact values after using the inequalities for $\La^{(s)}\pm\delta_{k+r^{(s)},i}\in\mathscr S_{\la^{(s)},\mu^{(s)}}$ are
\beq\label{eq:skew_orders_table}
\begin{array}{c|cc}
&\text{first factor}&\text{second factor}\\ \hline
k>m,\ i\le m+r^{(s)}&-1&0\\
k>m,\ i>m+r^{(s)}&0&-1.
\end{array}
\eeq

For example, for $i\le m+r^{(s)}$ and $\La^{(s)}+\delta_{k+r^{(s)},i} \in\mathscr S_{\la^{(s)},\mu^{(s)}}$,$\lambda_{k+r^{(s)}+1,i}=\lambda_{k+r^{(s)},i}+1$ and $\lambda_{k+r^{(s)}-1,i}
=\lambda_{k+r^{(s)},i}$, so the $j=i$ terms in \eqref{eq:skew_order_even} contribute $-1$; the terms with $j\ne i$ vanish by the strict order of the shifted even entries and the membership condition for the change.

To check the second factor put $c=l_{k+r^{(s)},i}^{(s)}$ and $p_0=m+r^{(s)}+l_{k+r^{(s)},i}^{(s)}$.
Since
\[
p_0\ge\lambda_{k+r^{(s)},m+r^{(s)}}+1,\quad \lambda_{k+r^{(s)}+1,m+r^{(s)}}\le p_0,
\]
\eqref{eq:skew_order_odd} equal to
\[
\begin{aligned}
&\de_{k+r^{(s)}-(m+r^{(s)})\ge p_0}+\de_{k+r^{(s)}-1-(m+r^{(s)})\ge p_0}\\
&\quad-\de_{k+r^{(s)}-(m+r^{(s)})\ge p_0+1}-\de_{k+r^{(s)}-(m+r^{(s)})\ge p_0}=0.
\end{aligned}
\]

For $\La^{(s)}-\delta_{k+r^{(s)},i}
\in\mathscr S_{\la^{(s)},\mu^{(s)}}$
one uses
$\lambda_{k+r^{(s)}+1,i}
=\lambda_{k+r^{(s)},i}$ and
$\lambda_{k+r^{(s)}-1,i}
=\lambda_{k+r^{(s)},i}-1$
and obtains the same values in~\eqref{eq:skew_orders_table}.
For $i>m+r^{(s)}$ the hook inequality
$\lambda_{k+r^{(s)},m+r^{(s)}}
\ge i-(m+r^{(s)})$,
combined with
$\lambda_{k+r^{(s)},i}\ge0$,
shows that every linear factor in the first product
of~\eqref{eq:skew_phi_odd} is nonzero at $z=-c$;
the interlacing of the odd entries leaves exactly
one uncancelled denominator in the second product.
Thus for $\La^{(s)}\pm\delta_{k+r^{(s)},i}\in\mathscr S_{\la^{(s)},\mu^{(s)}}$, $R_{k,\La^{(s)}}(z)$ has a simple pole at $z=-c^\pm_{k,s,i}$ with nonzero residue.

For clarity, if the numerator and denominator in any of \eqref{eq:skew_phi_even}--\eqref{eq:skew_phi_odd} are written as $\prod_\alpha(z+\alpha)$
and $\prod_\beta(z+\beta)$, with multiplicities, then we have
\beq\label{eq:skew_residue_product}
\operatorname*{Res}_{z=-c}R_{k,\La^{(s)}}(z)
=\frac{\displaystyle\prod_{\alpha\ne c}(\alpha-c)}
{\displaystyle\prod_{\beta\ne c}(\beta-c)}
\ne0.
\eeq

This removes all common factors before evaluation.
In particular, no indeterminate $0/0$ is used at an odd boundary.

The possible values of $c$ in~\eqref{eq:skew_pole_cases} for distinct changes of a fixed row which preserve $\mathscr S_{\la^{(s)},\mu^{(s)}}$ are different.

For positions $i,j\le m+r^{(s)}$ this follows from
$l_{k+r^{(s)},1}>
\cdots>l_{k+r^{(s)},m+r^{(s)}}$;
for $i,j>m+r^{(s)}$ from
$l_{k+r^{(s)},m+r^{(s)}+1}
<\cdots<l_{k+r^{(s)},k+r^{(s)}}$.
A coincidence between an increase and a decrease
in adjacent positions would force an equality
in the same row forbidden by the membership condition
for that change.
For $i\le m+r^{(s)}<j$, we have
\beq\label{eq:skew_no_collision}
\begin{aligned}
&l_{k+r^{(s)},i}-l_{k+r^{(s)},j}\\
&=\lambda_{k+r^{(s)},i}
+\lambda_{k+r^{(s)},j}
+(m+r^{(s)})-i
-\bigl(j-(m+r^{(s)})\bigr)+1,
\end{aligned}
\eeq
which, with the relevant increase or decrease inequalities,
excludes the remaining coincidences.
Different tensor factors have poles in distinct
cosets $-h^{(s)}+\BZ$.

Let $\Omega=\La+\delta_{k+r^{(s)},i}$
and $a=\ga_k-h^{(s)}-r^{(s)}-c^+_{k,s,i}$.
Taking the $\xi_\La$-diagonal coefficient
in $[e_k(v),f_k^{(1)}]$, and then the residue
at $v=a$, gives
\beq\label{eq:skew_product_plus}
C^+_{k;\Omega,\La}C^-_{k;\La,\Omega}
=(-1)^{\bk}
\operatorname*{Res}_{v=a}\Phi_{k,\La}(v).
\eeq

For $\Omega=\La-\delta_{k+r^{(s)},i}$,
with $a=\ga_k-h^{(s)}-r^{(s)}-c^-_{k,s,i}$,
it gives
\beq\label{eq:skew_product_minus}
C^-_{k;\Omega,\La}C^+_{k;\La,\Omega}
=-(-1)^{\overline{k+1}}
\operatorname*{Res}_{v=a}\Phi_{k,\La}(v).
\eeq

By~\eqref{eq:skew_phi_product},
\beq\label{eq:skew_product_factorization}
\begin{aligned}
\operatorname*{Res}_{v=a}\Phi_{k,\La}(v)
={}&
\operatorname*{Res}_{z=-c}R_{k,\La^{(s)}}(z)\\
&\times\prod_{t\ne s}
R_{k,\La^{(t)}}\bigl(
a+h^{(t)}+r^{(t)}-\ga_k
\bigr).
\end{aligned}
\eeq

The first factor is the explicit nonzero product
\eqref{eq:skew_residue_product}; the remaining factors
are nonzero and finite since
$h^{(s)}-h^{(t)}\notin\BZ$
and all their zeros and poles
have integral differences from $-h^{(t)}$.
Consequently
\[
C^+_{k;\Omega,\La}C^-_{k;\La,\Omega}\ne0
\quad\text{and}\quad
C^-_{k;\Omega,\La}C^+_{k;\La,\Omega}\ne0
\]
whenever the corresponding changed pattern belongs to
$\mathscr S_{\la^{(s)},\mu^{(s)}}$.
\epf

\bth\label{thm:ras}
Suppose $\Omega$ differs from $\La$ only by
$\Omega^{(s)}=\La^{(s)}\pm\delta_{k+r^{(s)},i}
\in\mathscr S_{\la^{(s)},\mu^{(s)}}$. Then
\beq\label{eq:skew_one_step}
\xi_\Omega=\frac{P_\Omega x_k^{(1)}\xi_\La}
{C^+_{k;\Omega,\La}}
\quad\text{or}\quad
\xi_\Omega=\frac{P_\Omega y_k^{(1)}\xi_\La}
{C^-_{k;\Omega,\La}},
\eeq
according as the change is an increase or a decrease.
In both cases the denominator is nonzero.
\eth
\bpf
By~\eqref{eq:skew_raise_coefficient},
\eqref{eq:skew_lower_coefficient} and Theorem~\ref{thm:ksi},
$P_\Omega x_k^{(1)}\xi_\La=C^+_{k;\Omega,\La}\xi_\Omega$
and similarly for $y_k^{(1)}$.
\epf

\bth\label{thm:GT basis}
The vectors~\eqref{eigvec} form a Gelfand--Tsetlin basis of
$W$. The actions of $d_k(u)$, $x_k(u)$ and $y_k(u)$ are
\eqref{eq:skew_eigen_formula} and~\eqref{eq:skew_x_y_action},
with the nonzero coefficient products
\eqref{eq:skew_product_plus}--\eqref{eq:skew_product_factorization}.
\eth
\bpf
This follows from Proposition~\ref{thm:bas} and
Theorems~\ref{thm:ksi} and~\ref{thm:ras}.
\epf

\bco\label{irr thm}
If $h^{(s)}-h^{(t)}\notin\BZ$ for $s\ne t$, then
$W$ is irreducible.
\eco
\bpf
Let $0\ne U\subseteq W$ be a submodule. Since the joint
Gelfand--Tsetlin spectrum is simple, for some $\La^{(s)}\in\mathscr S_{\la^{(s)},\mu^{(s)}}$
$(1\le s\le N)$,
\[
0\ne P_\La U=\BC\xi_\La\subset U.
\]
For $\Omega$ differing from $\La$ in one entry and
$\Omega^{(s)}\in\mathscr S_{\la^{(s)},\mu^{(s)}}$,
Theorem~\ref{thm:ras}
gives the formula
\[
\xi_\Omega
=\frac{P_\Omega x_k^{(1)}\xi_\La}{C^+_{k;\Omega,\La}}
\quad\text{or}\quad
\xi_\Omega
=\frac{P_\Omega y_k^{(1)}\xi_\La}{C^-_{k;\Omega,\La}}.
\]
By the defining conditions of Gelfand--Tsetlin patterns, any two elements of $\mathscr S_{\lambda,\mu}$ can be connected by successive changes of single entries, with all intermediate patterns remaining in $\mathscr S_{\lambda,\mu}$. Applying this to each tensor factor, we can obtain any $\xi_\Omega$ from $\xi_\Lambda$ by successive applications of the raising and lowering operators. Hence $U=W$.
\epf
\section{Classification theorem}
Let $L(\la(u))$ be the simple highest weight module of $\Ymn$, i.e. it is generated by a highest weight vector satisfying that
\begin{enumerate}
\item $t_{ij}(u)v=0$ for $1\leq i<j\leq m+n$,
\item $t_{ii}(u)v=\la_i(u)v$ for $1\leq i\leq m+n$.
\end{enumerate}

\bth[\cite{Zhang2}]
$L(\la(u))$ is finite dimensional if and only if its highest weight $\la(u)$ satisfies that
\[
\frac{\la_k(u)}{\la_{k+1}(u)}=\frac{P_k(u+(-1)^{\bk})}{P_k(u)},
\]
\[
\frac{\la_m(u)}{\la_{m+1}(u)}=\frac{Q_0(u)}{Q_1(u)},
\]
for $1\leq k<m+n$ and $k\neq m$, where $P_k(u)$ is a polynomial in $u$ and
\[
Q_0(u)=\prod_{s=1}^{N}(u+a_s),\quad Q_1(u)=\prod_{s=1}^{N}(u-b_s),
\]
such that $Q_0(u)$ and $Q_1(u)$ are coprime.
\eth
In fact, we know that up to an automorphism $\omega_f$ of $\Ymn$, Drinfeld polynomials uniquely determine a finite dimensional simple $\Ymn$-module.
Now let us recall a theorem from \cite{NT2}.
Let $n=0$ in $\Ymn$ and $V$ be any irreducible finite dimensional module over $\Y(\gl_m)$.
Let $P_1(u),\ldots,P_{m-1}(u)$ be the Drinfeld polynomials corresponding to $V$.
For each $1\leq k\leq m-1$ consider the collection of zeros of the polynomial $P_k(-u)$
\[
\{z_{ki}|i=1,\ldots, \text{deg} P_k\}
\]
\bth\cite{NT2}\label{thm:NT}
The action in module $V$ of the subalgebra generated by the coefficients of $d_1(u),\ldots,d_m(u)$ is semisimple
if and only if for all $k_1\geq k_2$ we have $z_{k_1i_1}-z_{k_2i_2}\neq 0,1,\ldots, k_1-k_2$ unless $(k_1,i_1)=(k_2,i_2)$\footnote{Notice that the Drinfeld polynomial $P_i(u)$ in \cite{NT2} is $P_i(u-i+1)$ in this paper, so the condition may be a little different from the original one in \cite{NT2}}.
\eth

For two integers $a>b$, we denote the set $\{b,b+1,b+2,\ldots,a\}$ by $[b;a]$.
If the simple $\Ymn$-module $V$ is of the form
\[
V=L_{h^{(1)}}(\la^{(1)})\otimes L_{h^{(2)}}(\la^{(2)})\otimes\cdots\otimes L_{h^{(M)}}(\la^{(M)}).
\]
The Drinfeld polynomials of $V$ can be written as
\[
P_i(u)=\prod_{s=1}^{M}(u+\la_i^{(s)}-1+h^{(s)})(u+\la_i^{(s)}-2+h^{(s)})\cdots (u+\la_{i+1}^{(s)}+h^{(s)}),
\]
\[
P_{m+j}(u)=\prod_{s=1}^{M}(u-\la_{m+j+1}^{(s)}+h^{(s)})(u-\la_{m+j+1}^{(s)}-1+h^{(s)})\cdots(u-\la_{m+j}^{(s)}+1+h^{(s)}),
\]
\[
Q_0(u)=\prod_{s=1}^{M}(u+\la_m^{(s)}+h^{(s)}),
\quad
Q_1(u)=\prod_{s=1}^{M}(u-\la_{m+1}^{(s)}+h^{(s)}),
\]
for $1\leq i\leq m-1$ and $1\leq j\leq n-1$.

It follows from Theorem \ref{thm:NT}
that if $z_{k_1i_1}-h^{(s_1)}\in [\la_{k_1+1}^{(s_1)};\la_{k_1}^{(s_1)}-1]$ and $z_{k_2i_2}-h^{(s_2)}\in [\la_{k_2+1}^{(s_2)};\la_{k_2}^{(s_2)}-1]$ with $1\leq s_1\neq s_2\leq M$, then $z_{k_1i_1}-z_{k_2i_2}\neq 0,1,\ldots,k_1-k_2, k_1-k_2$
for all $1\leq k_2\leq k_1\leq m-1$.

For a positive integer $M$, we denote an element $(h^{(1)},\ldots,h^{(M)})$ in $\BC^M$ with $h^{(a)}-h^{(b)}\in\BZ$ for all $1\leq a,b\leq M$ by $\underline{h}$.
For a collection of covariant $\gl_{m|n}$-weights $\underline{\la}=(\la^{(1)},\ldots,\la^{(M)})$, with every $\la^{(s)}$ is non-trivial for $1\leq s\leq M$, let us consider $V_{\underline{h}}(\underline{\la})$ as a tensor product of evaluation module with covariant weight, i.e.,
\[
V_{\underline{h}}(\underline{\la})=L_{h^{(1)}}(\la^{(1)})\otimes\cdots\otimes L_{h^{(M)}}(\la^{(M)}).
\]
We first treat two nontrivial factors:
\[
 V=L_{h^{(1)}}(\la^{(1)})\otimes L_{h^{(2)}}(\la^{(2)}).
\]
When their parameter difference is integral, comparisons of shifted
weights below are comparisons of real integers after subtracting their
common complex part. A relabelling in an inequality does not permute
the actual tensor factors.

Recall the Gelfand-Tsetlin basis given in Section \ref{sec:GT bas}.
The basis of $L_{h^{(s)}}(\la^{(s)})$ can be parameterized by
the patterns $\La^{(s)}\in \mathscr{S}_{\la^{(s)}}$ in the form
\begin{equation*}
\resizebox{.9\textwidth}{!}{$\La^{(s)}=\begin{array}{cccccccc}
   \la_1^{(s)}   &  \cdots      & \la_m^{(s)} & \la_{m+1}^{(s)} & \cdots & \la_{m+n-1}^{(s)} & \la_{m+n}^{(s)} \\
  \la_{m+n-1,1}^{(s)} &  \cdots     & \la_{m+n-1,m}^{(s)} & \la_{m+n-1,m+1}^{(s)} & \cdots & \la_{m+n-1,m+n-1}^{(s)} & \\
  \vdots    &  \vdots & \vdots & \vdots   & \reflectbox{$\ddots$}    \\
  \la_{m+1,1}^{(s)} &  \cdots  & \la_{m+1,m}^{(s)} & \la_{m+1,m+1}^{(s)}    \\
  \la_{m,1}^{(s)} &  \cdots  & \la_{m,m}^{(s)}                           \\
  \la_{m-1,1}^{(s)} &  \cdots                             \\
  \vdots  &  \reflectbox{$\ddots$}                              \\
  \la_{11}^{(s)}\\
\end{array}$}.
\end{equation*}
\medskip\noindent\textit{The even-side necessary condition.}
We first justify the restriction argument.
For $1\leq q\leq m+n$, Set $H_q=\sum_{i=1}^q E_{ii}$ and
\[
L_{h^{(s)}}(\la^{(s)})=\bigoplus_a (L_{h^{(s)}}(\la^{(s)}))_a,
\]
where $(L_{h^{(s)}}(\la^{(s)}))_a=\{v:H_qv=av\}$.
Consider the filtration of $V$ induced by the lexicographic order on the pairs
$(a,b)$, with larger weights first.
For $1\le i,j\le q$, we have
\beql{cop sum}
\Delta(t_{ij}(u))=\sum_{k=1}^{q}t_{ik}(u)\otimes t_{kj}(u)+\sum_{k=q+1}^{m+n}t_{ik}(u)\otimes t_{kj}(u),
\eeq
Due to
\[
[H_q,t_{ia}(u)]=t_{ia}(u),\quad [H_q,t_{aj}(u)]=-t_{aj}(u),
\]
for $i,j\leq q<a$.
Thus the second summand in \eqref{cop sum} maps
\[
V_a^{(1)}\otimes V_b^{(2)}\longrightarrow V_{a+1}^{(1)}\otimes V_{b-1}^{(2)}.
\]
Consequently, the induced action on the associated graded module is given by
\[
\operatorname{gr}\Delta(t_{ij}(u))=\sum_{a=1}^{q}t_{ia}(u)\otimes t_{aj}(u),\qquad 1\le i,j\le q,
\]
which is the standard coproduct of the corresponding (super-)Yangian.

Since the action of the Gelfand--Tsetlin subalgebra is semisimple on $V$, so is $\operatorname{gr}V$
for the corresponding leading Gelfand--Tsetlin subalgebra. Hence every simple subquotient of a
tensor product of constituents selected from the
two restricted factors is tame.

For $q=m$, the possible highest weights $\alpha^{(s)}$ are precisely the
partitions satisfying
\beql{legal-even-row}
 \max(0,\la_i^{(s)}-n)\le\alpha_i^{(s)}\le\la_i^{(s)},
 \qquad \alpha_1^{(s)}\ge\cdots\ge\alpha_m^{(s)}\ge0.
\eeq
Equivalently these are the possible rows $\la_{m,i}^{(s)}$ in
$\mathscr S_{\la^{(s)}}$. The tensor of their highest vectors in the graded
restriction has a simple highest weight quotient with
\beql{even-branch-polynomials}
 P'_i(u)=\prod_{s=1}^{2}
       \prod_{a=\alpha_{i+1}^{(s)}}^{\alpha_i^{(s)}-1}(u+h^{(s)}+a).
\eeq
Theorem~\ref{thm:NT} applies to this quotient.

For $n=1$, assume that $\lambda_1^{(s)}+h^{(s)}\geq\lambda_1^{(t)}+h^{(t)}$ $(s,t=1,2)$,
and put $k_s=\max\{i\leq m:\lambda_i^{(s)}>0\}$.
Applying Theorem~5.2 to \eqref{legal-even-row}, with the $m$-th rows varying over
$\mathscr S_{\lambda^{(s)}}$, gives
\beql{weak-even-bound}
 \max(0,\la_m^{(s)}-1)+h^{(s)}-\la_1^{(t)}-h^{(t)}
 \ge\begin{cases}k_s-1,&k_s<m,\\m-2,&k_s=m.\end{cases}
\eeq

Suppose now that $k_s=m$ and equality with $m-2$ holds in
\eqref{weak-even-bound}. Write $\la=\la^{(s)}$, $\nu=\la^{(t)}$,
$c=\la_{m+1}$, $d=\nu_{m+1}$.
Thus
\beql{mixed-endpoint}
 h^{(s)}+\la_m=h^{(t)}+\nu_1+m-1.
\eeq

Choose the following \((m-1)\)-st rows
\[
\mu^+=(\lambda_1,\ldots,\lambda_{m-1}),\qquad
\mu^-=(\max\{\nu_2-1,0\},\ldots,\max\{\nu_m-1,0\}),
\]
and complete the lower rows of the corresponding Gelfand--Tsetlin patterns by taking their maximal permitted entries. The latter extends through the \(m\)-th row
\[
(\nu_1,\max\{\nu_2-1,0\},\ldots,\max\{\nu_m-1,0\}).
\]
The shifted entries of the first factor are greater than \(h^{(t)}+\nu_1\), since
\[
h^{(s)}+\lambda_m-m+2=h^{(t)}+\nu_1+1.
\]
Hence the joint eigenvalues of $B_1(u),\ldots,B_{m-1}(u)$ distinguish the chosen Gelfand--Tsetlin patterns through level $m-1$. Under the tame assumption, the corresponding joint eigenspace is nonzero and invariant under $\psi_{m-1}(\Y(\mathfrak{gl}_{1|1}))$.

For one factor with $(m-1)$-st row $\mu$, the ratio of the highest-weight eigenvalues of $\psi_{m-1}(t_{11}(u))$ and $\psi_{m-1}(t_{22}(u))$ is
\[
\frac{\displaystyle\prod_{i=1}^{m}(u+h+\lambda_i+m-i)}
{(u+h-\lambda_{m+1})
\displaystyle\prod_{i=1}^{m-1}(u+h+\mu_i+m-i)}.
\]
For the two chosen $(m-1)$-th row these ratios reduce to
\[
 \frac{u+h^{(s)}+\lambda_m}{u+h^{(s)}-c},\qquad
 \begin{cases}
 \displaystyle\frac{u+h^{(t)}+\nu_1+m-1}{u+h^{(t)}-d},&k_t=m,\\[5pt]
 \displaystyle\frac{u+h^{(t)}+\nu_1+m-1}{u+h^{(t)}+m-k_t},&k_t<m.
 \end{cases}
\]
The highest-weight ratio of the corresponding simple quotient is the product of the two expressions above. By \eqref{mixed-endpoint} the numerator
has a double root. Neither denominator cancels it, since
\[
\lambda_m+c>0,\qquad \nu_1+d+m-1>0\ (k_t=m),\qquad \nu_1+k_t-1>0\ (k_t<m).
\]
But a simple $\Y(\gl_{1|1})$-module is tame only if the numerator of its reduced highest weight ratio has simple roots \cite[Corollary 3.9 and Proposition 3.10]{Lu}. This is a contradiction.
For $m=1$ the same argument applies directly to the two highest vectors: failure of the required inequality is precisely equality of their numerator roots. We have proved the stronger bound with $k_s-1$ also when $k_s=m$.

For general $n\ge1$ and $m\ge2$, the same even-row argument uses
\[
 h^{(s)}+[\max(0,\lambda_{i+1}-n);\lambda_i-1]
\]
in place of the preceding intervals. Their union has lower endpoint
$h^{(s)}$ if $k_s<m$, and $h^{(s)}+\max(0,\lambda_m-n)$ otherwise.
Comparison with the largest root at index $1$ of the other factor gives
\[
 D:=\max(0,\lambda_m-n)+h^{(s)}-\nu_1-h^{(t)}
 \ge\begin{cases}k_s-1,&k_s<m,\\m-2,&k_s=m.\end{cases}
\]
Only $k_s=m$, $D=m-2$ remains. Use the leading
$\Y(\gl_{m|1})$ filtration described above. In the first factor choose the covariant constituent with even entries
\[
 \beta_i=\begin{cases}
 \max(1,\lambda_i-n+1),&\lambda_i>0,\\0,&\lambda_i=0,
 \end{cases}
 \qquad\beta_{m+1}=\lambda_{m+n}.
\]
It belongs to the branching set: for $1\le j\le n$ it extends through the rows
\[
 \lambda_{m+j,i}=\begin{cases}
 \max(1,\lambda_i-n+j),&\lambda_i>0,\\0,&\lambda_i=0,
 \end{cases}
 \quad
 \lambda_{m+j,m+a}=\lambda_{m+n-j+a}\quad(1\le a\le j).
\]
For the second factor choose the constituent
$(\nu_1,\ldots,\nu_m\mid\nu_{m+1})$, extending it by the first $j$ odd
entries at row $m+j$. At this remaining endpoint, putting $L=\max(0,\lambda_m-n)$, one has
\[
 \beta_m=L+1,\qquad
 \beta_m+h^{(s)}=\nu_1+h^{(t)}+m-1.
\]
The order required for the $n=1$ endpoint argument is therefore the
correct one, and that argument excludes the endpoint.

If $m=1$, failure of the desired inequality means that
$\max(1,\lambda_1-n+1)\le \nu_1+h^{(t)}-h^{(s)}\le\lambda_1$.
Choose the following leading $\gl_{1|1}$ constituent
$(\nu_1+h^{(t)}-h^{(s)}|\lambda_{1+n})$ in factor $s$, and
$(\nu_1\mid\nu_2)$ in factor $t$. The first extends through even entries
$\max(\nu_1+h^{(t)}-h^{(s)},\lambda_1-n+j)$ and the last $j$ odd entries of $\lambda$.
Their numerator roots coincide, and neither denominator cancels this
root because both even highest entries are positive. The rank-one
criterion again gives a contradiction. We have thus proved, for all
$m,n\ge1$,
\beql{eq:even_inequality}
\max(0,\la_m^{(s)}-n)+h^{(s)}-\la_1^{(t)}-h^{(t)}\ge k_s-1\quad \text{ if } \lambda_1^{(s)}+h^{(s)}\geq\lambda_1^{(t)}+h^{(t)}.
\eeq

\medskip\noindent\textit{The necessary conditions on the odd side.}
Set:
\beql{row-integers}
 \nu_i^{(s)}=\begin{cases}
 \la_i^{(s)},&1\le i\le m,\\
 \#\{1\le j\le n:\la_{m+j}^{(s)}\ge i-m\},&i>m,
 \end{cases}
\eeq
for $s=1,2$.
For $0\le r\le n-2$ consider
indices $i$ with $\nu_i^{(s)}>0$ and integers $a^{(s)}$ satisfying
\beql{eq:odd_optional}
\max(0,\nu_i^{(s)}-n+r+1)\le a^{(s)}\le
\begin{cases}
\nu_i^{(s)}-1,&i\le m,\\
\min(r,\nu_i^{(s)}-1),&i>m.
\end{cases}
\eeq
The additional condition is the following comparison of shifted entries:
\beql{eq:odd_shift_comparison}
 h^{(1)}+a^{(1)}-i\ne h^{(2)}+a^{(2)}-j,
\eeq
for every $r$ and every pair $(i,a^{(1)}),(j,a^{(2)})$ in the ranges
\eqref{eq:odd_optional} of their respective factors.

We verify both pattern membership and necessity. A row at level $m+r$ specifies
a nonincreasing sequence $\alpha_i^{(s)}$ by
\[
 \alpha_i^{(s)}=\begin{cases}
 \la_{m+r,i}^{(s)},&i\le m,\\
 \#\{1\le a\le r:\la_{m+r,m+a}^{(s)}\ge i-m\},&i>m.
 \end{cases}
\]
Then we have the following inequality.
\beql{legal-odd-prefix}
 0\le\nu_i^{(s)}-\alpha_i^{(s)}\le n-r,
 \qquad\alpha_{m+1}^{(s)}\le r,
 \qquad\alpha_i^{(s)}\ge\alpha_{i+1}^{(s)}\ge0.
\eeq
The entries between the $m+r$-th row and $m+n$-th row are determined when $\nu_i^{(s)}-\alpha_i^{(s)}=0$ or $n-r$.
Otherwise, $1\leq\nu_i^{(s)}-\alpha_i^{(s)}\leq n-r-1$, which gives the range of $a^{(s)}=\alpha_i^{(s)}$
in \eqref{eq:odd_optional}.

Assume \eqref{eq:even_inequality}, and let $r$ be the first level at which \eqref{eq:odd_shift_comparison} fails. By the preceding separation conditions, the joint eigenvalues of $B_1(u),\ldots,B_{m+r}(u)$ determine the corresponding rows of patterns in each factor. Choose such rows realizing the failure of \eqref{eq:odd_shift_comparison}.

Since $V$ is tame, their joint eigenspace is nonzero and invariant under the subalgebra generated by $d_{m+r+j}(u)$ and the corresponding raising and lowering operators. Replacing $u$ by $-v$ identifies the resulting all-odd action with the $\Y(\mathfrak{gl}_{n-r})$-action. The Drinfeld polynomials of its simple highest-weight quotient are
\beql{eq:tail_Drinfeld}
P'_j(v)=\prod_{s=1}^{2}\prod_{\nu_i^{(s)}-\alpha_i^{(s)}=j}\bigl(v+i-\alpha_i^{(s)}-m+r-1-h^{(s)}\bigr),\qquad1\le j<n-r.
\eeq
Indeed, an index satisfying
$\nu_i^{(s)}-\alpha_i^{(s)}=j$
contributes the factor
\[
\frac{v+i-\alpha_i^{(s)}-m+r-1-h^{(s)}+1}{v+i-\alpha_i^{(s)}-m+r-1-h^{(s)}},
\]
to the first $j$ diagonal components.
Hence the failure of \eqref{eq:odd_shift_comparison} yields either a repeated root or a common root
of two polynomials in \eqref{eq:tail_Drinfeld}, contradicting Theorem \ref{thm:NT}.
This proves the necessity.

\bth\label{thm:ns con}
Suppose that $m,n> 0$ , $\la^{(1)},\la^{(2)}$ are nontrivial covariant weights of $\gl_{m|n}$ and
 
\[
V=L_{h^{(1)}}(\la^{(1)})\otimes L_{h^{(2)}}(\la^{(2)})
\]
is a simple  $\Y(\mathfrak{gl}_{m|n})$-module. Then
$V$ is a tame $\Y(\mathfrak{gl}_{m|n})$-module if and only if $h^{(1)}-h^{(2)}\notin\BZ$ or $h^{(1)}-h^{(2)}\in\BZ$ and the following conditions hold: for distinct $s,t\in \{1,2\}$, ordered so that
\[
\la_1^{(s)}+h^{(s)}\geq\la_1^{(t)}+h^{(t)},
\]
we have
\[
\max\{0,\la_m^{(s)}-n\}+h^{(s)}-\la_1^{(t)}-h^{(t)}\geq k_s-1,
\]
where
\[
k_s=\max\{i\leq m:\la_i^{(s)}>0\}.
\]
Moreover,
\[
h^{(s)}+a^{(s)}-i\neq h^{(t)}+a^{(t)}-j
\]
for all $(i,a^{(s)})$ and $(j,a^{(t)})$ in the respective ranges \eqref{eq:odd_optional}.
\eth
\bpf
The necessity was proved above. We prove the sufficiency by showing that the joint eigenvalues
of $B_k(u)$, $1\leq k\leq m+n$, distinguish all pairs of Gelfand--Tsetlin patterns.

For $1\leq k\leq m$, the diagonal eigenvalues are
\[
\prod_{s=1}^{2}\prod_{i=1}^{k}
(u+h^{(s)}+l_{ki}^{(s)}).
\]
If $h^{(1)}-h^{(2)}\in\BZ$ by \eqref{eq:even_inequality} and the interlacing conditions of Section \ref{sec:GT bas},
then for $1\leq i\leq\min\{k,k_s\}$,
\[
\begin{aligned}
l_{ki}^{(s)}+h^{(s)}
&\geq\lambda_{mm}^{(s)}-i+1+h^{(s)}\\
&\geq\max\{0,\lambda_m^{(s)}-n\}
-i+1+h^{(s)}\\
&\geq\lambda_1^{(t)}+h^{(t)}+k_s-i\\
&\geq l_{k1}^{(t)}+h^{(t)}.
\end{aligned}
\]
For $k>k_s$, the entries $\lambda_{ki}^{(s)}=0$ $(i>k_s)$ are fixed by restriction on patterns.
After removing their known factors and accounting for the boundary equalities, the diagonal
eigenvalues determine the even rows of $\La^{(1)}$ and $\La^{(2)}$ separately.

Now suppose that the rows through level $m+r$ are determined. By Lemma \ref{lem:chkL}, each undetermined
increment in the next row contributes
\[
\frac{u+c+1}{u+c},\qquad
c=h^{(s)}+\alpha_i^{(s)}+1-i.
\]
Within each factor, these values of $c$ are distinct, and \eqref{eq:odd_shift_comparison} separates those belonging
to different factors. If two pairs of patterns
have the same joint eigenvalues, comparison
of the next rows gives
\[
\prod_{c\in F}
\left(\frac{u+c+1}{u+c}\right)^{e_c}=1,
\qquad e_c\in\BZ,\quad F\text{ finite}.
\]
Comparing the orders at $u=-c$ gives $e_c=e_{c-1}$. Since $F$ is finite, we obtain $e_c=0$ for all $c$.
Hence every increment in the next row is uniquely determined. At $r=n-1$ the remaining
increments are forced. Induction on $r$ proves that all joint eigenvalues are distinct.

Clearly, if $h^{(1)}-h^{(2)}\notin\BZ$, the above statement about that the joint eigenvalues are distinct automatically holds.

Consequently, by Proposition \ref{copro on B} and the spectral projection construction \eqref{eq:skew_A}-\eqref{eq:skew_leading},
$V$ admits a Gelfand--Tsetlin eigenbasis with simple spectrum. Hence $V$ is tame.
\epf

\bth\label{thm:tame_multiple}
Let $\la^{(s)}$ ($1\leq s\leq M$) be nontrivial covariant weights of $\gl_{m|n}$, suppose that
\[
V=L_{h^{(1)}}(\la^{(1)})\otimes\cdots\otimes
L_{h^{(M)}}(\la^{(M)})
\]
is a simple  $\Y(\mathfrak{gl}_{m|n})$-module.
Then the following conditions are equivalent:
\begin{enumerate}
\item $V$ is a tame $\Y(\mathfrak{gl}_{m|n})$-module.

\item For every $1\leq s<t\leq M$, the
tensor product
\[
L_{h^{(s)}}(\la^{(s)})
\otimes L_{h^{(t)}}(\la^{(t)})
\]
is a tame $\Y(\mathfrak{gl}_{m|n})$-module.
\end{enumerate}
Moreover, under these conditions, $V$ has simple Gelfand--Tsetlin spectrum.
\eth

\bpf
We first prove $(1)\Rightarrow(2)$.
Let $\sigma$ be any permutation of
$\{1,\ldots,M\}$, and put
\[
V_\sigma=L_{h^{(\sigma(1))}}(\la^{(\sigma(1))})\otimes\cdots\otimes
L_{h^{(\sigma(M))}}(\la^{(\sigma(M))}).
\]
The tensor product of the highest weight vectors generates a submodule of $V_\sigma$ whose simple
highest weight quotient is isomorphic to $V$.
Since $\dim V_\sigma=\dim V$, we obtain
\[
V_\sigma\cong V.
\]
Every consecutive subproduct of $V_\sigma$ is also irreducible, since a proper nonzero
submodule would give a proper nonzero submodule of $V_\sigma$ after tensoring with the remaining
factors.

We next show that tameness passes to consecutive subproducts. Let $U$ and $Z$ be nonzero
finite-dimensional weight modules such that $U\otimes Z$ is tame. Write
\[
U=\bigoplus_\mu U_\mu.
\]
Use the standard RTT normalization $d_i(u)=1+O(u^{-1})$ in the following calculation. The polynomial normalization differs by known scalar powers of $u$ and does not affect semisimplicity.

By the triangular Gauss coproduct in Proposition \ref{copro on B}, the filtration of $U\otimes Z$
induced by the weights of $U$ has
\[
\left.\operatorname{gr}d_i(u)\right|
_{U_\mu\otimes Z}
=
\left.d_i(u)\right|_{U_\mu}
\otimes\left.d_i(u)\right|_Z.
\]
Since $U\otimes Z$ is tame, the commuting coefficients of $d_i(u)$ act semisimply on every associated graded component.

Choose a common eigenvector $0\neq z\in Z$ such that
\[
d_i(u)z=\zeta_i(u)z.
\]
On $U_\mu\otimes\BC z$, the induced action is
\[
\zeta_i(u)\left.d_i(u)\right|_{U_\mu}.
\]
As $\zeta_i(u)$ is invertible, multiplication by this scalar series does not change the commutative algebra generated by the coefficients. Thus its action on $U_\mu$ is semisimple.
Since this holds for every $\mu$, the module $U$ is tame.
Choosing a common eigenvector in a nonzero $U_\mu$ proves similarly that $Z$ is tame.

It follows that every consecutive subproduct of $V_\sigma$ is tame. For any $s<t$, choose
$\sigma$ placing these two factors consecutively.
Then
\[
L_{h^{(s)}}(\la^{(s)})\otimes L_{h^{(t)}}(\la^{(t)})
\]
is simple and tame, proving (2).

Now, if
\[
L_{h^{(s)}}(\la^{(s)})\otimes L_{h^{(t)}}(\la^{(t)})
\]
is simple and tame, then the conditions shown in Theorem \ref{thm:ns con} hold.
Similar to the proof of Theorem \ref{thm:ns con},
We can show that these joint characters distinguish all tuples of Gelfand--Tsetlin patterns.

Consequently, distinct tuples of patterns have distinct joint $B_k(u)$-characters.
By Proposition \ref{copro on B} and the spectral projection construction \eqref{eq:skew_A}-\eqref{eq:skew_leading},
$V$ admits a Gelfand--Tsetlin eigenbasis with simple spectrum. Hence $V$ is tame.
\epf

Take the tensor product of the Gelfand-Tsetlin basis vectors of each $L_{h^{(s)}}(\la^{(s)})\ (1\leq s\leq M)$ as a basis of $V$.  In this basis, the matrix realization of $d_i(u)$ is given by
\[
\begin{pmatrix}
\alpha_{i1}&*&*&\cdots&*\\
0&\alpha_{i2}&*&\cdots&*\\
\vdots&\vdots&\vdots&&\vdots\\
0&0&0&\cdots&\alpha_{iK}
\end{pmatrix}.
\]
If for any $1\leq a\neq b\leq K$, there is $1\leq i\leq m+n$, such that $\alpha_{ia}\neq \alpha_{ib}$, we say that $(h^{(s)},\la^{(s)})$ $(1\leq s\leq M)$ satisfying the strong non-crossing condition.
\bre
If $n=0$, \cite{NT2} give a necessary and sufficient condition for a finite dimensional simple $\Y(\gl_m)$-module to be tame. We can also translate that to a condition by comparing the highest weights: for every pair of distinct indices $s,t$ and $1\leq l\leq k\leq m-1$ such that
\[
\la_l^{(s)}>\la_{l+1}^{(s)},\qquad \la_k^{(t)}>\la_{k+1}^{(t)},
\]
we have either
\[
h^{(s)}+\la_l^{(s)}+k-l\leq h^{(t)}+\la_{k+1}^{(t)}
\]
or
\[
h^{(t)}+\la_k^{(t)}\leq h^{(s)}+\la_{l+1}^{(s)}.
\]
If $m=0$, the corresponding conditions are obtained by replacing $m$ by $n$ and $h^{(s)}$ by $-h^{(s)}$.

In fact, either $mn=0$ or not, combining Theorem \ref{thm:ns con} and results from \cite{NT2}, the strong non-crossing condition is equivalent to the condition given above or the one shown in Theorem \ref{thm:ns con} (2).
Thus under the assumption in Theorem \ref{thm:tame_multiple},
the following two statements are equivalent:
\begin{enumerate}
\item[(a)] $d_i(u)$ ($1\leq i\leq m+n$) act semisimply on $V$,
\item[(b)] $(h^{(s)},\la^{(s)})$ $(1\leq s\leq M)$ satisfies the strong non-crossing condition.
\end{enumerate}
\ere

\section{Connection between $V_{\underline{h}}(\underline{\la})$ and skew modules}
In this section, we consider the connection between $V$ constructed in Theorem \ref{thm:tame_multiple} and the skew modules. For simplicity, we may set the index $t$ in Theorem \ref{thm:tame_multiple} to be $s+1$.

Let us consider the following module of $\Ymn$:
\[
L_{h}(\la\slash \mu)\otimes L_{h}(\la^{(1)})\otimes L_{h}(\la^{(2)})\otimes\cdots\otimes L_{h}(\la^{(M-1)}),
\]
such that for $1\leq b\leq M-1$,
\[
\la^{(b)}=(\la_m^{(b)},\la_m^{(b)},\ldots,\la_m^{(b)},\la_{m+1}^{(b)},\la_{m+2}^{(b)},\ldots,\la_{m+n}^{(b)})
\]
\[
\begin{split}
\la=&(\la_1,\la_2,\ldots,\la_{m-1},\underbrace{\la_m^{(1)}-p_{M-1},\ldots,\la_m^{(1)}-p_{M-1}}_{p_1-m+2},\\
&\la_{p_1+2},\ldots,\la_{p_1+m-1},\underbrace{\la_m^{(2)}-p_{M-1}+p_1,\ldots,\la_m^{(2)}-p_{M-1}+p_1}_{p_2-p_1-m+2},\\
&\ldots,\\
&\la_{p_{M-2}+2},\ldots,\la_{p_{M-2}+m-1},\underbrace{\la_m^{(M-1)}-p_{M-1}+p_{M-2},\ldots,\la_m^{(M-1)}-p_{M-1}+p_{M-2}}_{p_{M-1}-p_{M-2}-m+2},\\
&\la_{p_{M-1}+2},\ldots,\la_{p_{M-1}+m},\la_{p_{M-1}+m+1},\ldots,\la_{p_{M-1}+m+n})
\end{split}
\]
and
\[
\begin{split}
\mu=&(\underbrace{\la_m^{(1)}-p_{M-1},\ldots,\la_m^{(1)}-p_{M-1}}_{p_1},\underbrace{\la_m^{(2)}-p_{M-1}+p_1,\ldots,\la_m^{(2)}-p_{M-1}+p_1}_{p_2-p_1},\\
&\ldots,\underbrace{\la_m^{(M-1)}-p_{M-1}+p_{M-2},\ldots,\la_m^{(M-1)}-p_{M-1}+p_{M-2}}_{p_{M-1}-p_{M-2}}),
\end{split}
\]
for a sequence of integers $p_b\geq p_{b-1}+m-1$ with $p_0=0$. Here, $L_h(\la\slash \mu)$ is a skew module of $\Ymn$ through the homomorphism $\pi_{m+p_{M-1}|n}\psi_{p_{M-1}}$ defined in Section \ref{sec:skew mod}.

We denote such a collection $(\la\slash\mu,\la^{(1)},\la^{(2)},\cdots,\la^{(M-1)})$ by $\Theta(\underline{\la})$ and denote
\[
\mathcal{M}_h(\Theta(\underline{\la}))=L_h(\la\slash \mu)\otimes L_{h}(\la^{(1)})\otimes L_{h}(\la^{(2)})\otimes\cdots\otimes L_{h}(\la^{(M-1)}).
\]
Suppose that $\xi_0$, $\xi^{(s)}_0$ are highest weight vectors in $L_h(\la\slash \mu)$ and $L_h(\la^{(s)})$ ($1\leq s\leq M-1$), respectively.
Let $\mathcal{L}_h(\Theta(\underline{\la}))$ be the simple quotient of the $\Ymn$-module generated by $\xi_0\otimes \xi^{(1)}_0\otimes\cdots\otimes \xi^{(M-1)}_0$.
\bre
A $\gl_{m+p_{M-1}|n}$ weight $\la$ in the above form satisfies $\la_{m}^{(b)}-p_{M-1}+p_{b-1}\geq\la_{p_b+2}\geq \la_{p_b+m-1}\geq \la_m^{(b+1)}-p_{M-1}+p_b$.
\ere

\bth\label{thm:ske ten}
Let $\underline{\la}=(\la^{(1)},\ldots,\la^{(M)})$ is a collection of covariant weight and
\[
V_{h}(\underline{\la})=L_{h}(\la^{(1)})\otimes L_{h}(\la^{(2)})\otimes\cdots\otimes L_{h}(\la^{(M)}).
\]
be a simple and tame $\Ymn$-module.
Then $V_h(\underline{\la})$ coincides with $\mathcal{L}_h(\Theta(\underline{\la}))$ up to some automorphism $\omega_f$ of $\Ymn$.
\eth
\bpf

We already know that the Drinfeld polynomials of $V_{h}(\underline{\la})$ are
\[
P_i(u)=\prod_{s=1}^{M}(u+\la_i^{(s)}-1+h)(u+\la_i^{(s)}-2+h)\cdots (u+\la_{i+1}^{(s)}+h),
\]
\[
P_{m+j}(u)=\prod_{s=1}^{M}(u-\la_{m+j+1}^{(s)}+h)(u-\la_{m+j+1}^{(s)}-1+h)\cdots(u-\la_{m+j}^{(s)}+1+h),
\]
\[
Q_0(u)=\prod_{s=1}^{M}(u+\la_m^{(s)}+h),
\quad
Q_1(u)=\prod_{s=1}^{M}(u-\la_{m+1}^{(s)}+h),
\]
for $1\leq i\leq m-1$ and $1\leq j\leq n-1$.

Let us construct $\Theta(\underline{\la})$ such that Drinfeld polynomials of $\mathcal{L}_h(\Theta(\underline{\la}))$ are exactly  the ones above.

Due to Theorem \ref{thm:tame_multiple}, $P_i(u)$ has no multiple roots. Without loss of generality, we assume that $\la_{i+1}^{(s)}>\la_i^{(s+1)}-1$ for $1\leq s\leq M-1$. Furthermore, we iteratively define a sequence $q_1,q_2,\ldots,q_M$ of integers by setting  $q_1=1$ and
\[
\la_m^{(s)}-\la_1^{(s+1)}=q_{s+1}-q_s,
\]
for $1\leq s\leq M-1$. Again by Theorem \ref{thm:tame_multiple}, $q_{s+1}-q_s\geq m-1$.
Set $r=q_M-1$.

Take a $\gl_{m+r|n}$-weight $\la$ by setting
\[
\begin{split}
\la=&(\la_1^{(1)}-q_M+q_1,\la_2^{(1)}-q_M+q_1,\ldots,\la_{m-1}^{(1)}-q_M+q_1,\underbrace{\la_m^{(1)}-q_M+q_1,\ldots,\la_m^{(1)}-q_M+q_1}_{q_2-q_1-m+2},\\
&\la_2^{(2)}-q_M+q_2,\ldots,\la_{m-1}^{(2)}-q_M+q_2,\underbrace{\la_m^{(2)}-q_M+q_2,\ldots,\la_m^{(2)}-q_M+q_2}_{q_3-q_2-m+2},\\
&\vdots\\
&\la_2^{(M-1)}-q_M+q_{M-1},\ldots,\la_{m-1}^{(M-1)}-q_M+q_{M-1},\underbrace{\la_m^{(M-1)}-q_M+q_{M-1},\ldots,\la_m^{(M-1)}-q_M+q_{M-1}}_{q_{M}-q_{M-1}-m+2},\\
&\la_2^{(M)},\ldots,\la_{m}^{(M)},\la_{m+1}^{(M)},\la_{m+2}^{(M)},\ldots,\la_{m+n}^{(M)}).
\end{split}
\]
Then take a $\gl_r$-weight $\mu$ by setting
\[
\begin{split}
\mu=&(\underbrace{\la_m^{(1)}-q_M+q_1,\ldots,\la_m^{(1)}-q_M+q_1}_{q_2-q_1},\underbrace{\la_m^{(2)}-q_M+q_2,\ldots,\la_m^{(2)}-q_M+q_2}_{q_3-q_2},\ldots,\\
& \underbrace{\la_m^{(M-1)}-q_M+q_{M-1},\ldots,\la_m^{(M-1)}-q_M+q_{M-1}}_{q_{M}-q_{M-1}}).
\end{split}
\]
As a result, we get a skew module $L_h(\la\slash\mu)$ of $\Ymn$ as defined in Section \ref{sec:skew mod}.

Set $\Theta(\underline{\la})=(\la\slash\mu,\la^{(1)},\ldots,\la^{(M-1)})$. Then $\mathcal{L}_h(\Theta(\underline{\la}))$ possesses the desired Drinfeld polynomials, i.e., $V_h(\underline{\la})$ is isomorphic to $\mathcal{L}_h(\Theta(\underline{\la}))$  up to some automorphism $\omega_f$ of $\Ymn$.
\epf
\bre
Let
\[
V=V_{h^{(1)}}(\underline{\la}^{(1)})\otimes\cdots\otimes V_{h^{(N)}}(\underline{\la}^{(N)})
\]
be the $\Ymn$-module for some $h^{(a)}\in\BC$ and $1\leq a\leq N$. Then, by Theorem \ref{thm:ske ten}, we may write $V$ as
\[
\mathcal{L}_{h^{(1)}}(\Theta(\underline{\la}^{(1)}))\otimes \mathcal{L}_{h^{(2)}}(\Theta(\underline{\la}^{(2)}))\otimes\cdots\otimes \mathcal{L}_{h^{(N)}}(\Theta(\underline{\la}^{(N)})).
\]

Notice that when $n=0$, the tensor factor $L_h(\la^{(a)})$ of $\mathcal{M}_h(\Theta(\underline{\la}))$ is  $1$-dimensional. Thus $\mathcal{M}_h(\Theta(\underline{\la}))$ degenerates into a skew module from \cite{NT2}. In this case, $\mathcal{M}_h(\Theta(\underline{\la}))\cong \mathcal{L}_h(\Theta(\underline{\la}))$ and our result coincides with Theorem 4.1 in \cite{NT2}.
\ere

\bigskip
\centerline{\bf Acknowledgments}
\medskip
V. Futorny is partially supported by  the National Natural Science Foundation of China (Grants 12350710787 and 12350710178).
J. Zhang is partially supported by the National Natural Science Foundation of China (Grant 12571026).
Z. Li is partially supported by the National Natural Science Foundation of China (Grant 12601052).


\begin{thebibliography}{99}
\bibitem{BR}
{A. Berele and A. Regev},
{\it Hook Young diagrams with applications to combinatorics and to representations of Lie superalgebras},
Adv. in Math. {\bf 64} (1987), no.~2, 118--175.

\bibitem{Cher}
{I. Cherednik},
{\it A new interpretation of Gel'fand-Tzetlin bases},
Duke Math. J. {\bf 54} (1987), no.~2, 563--577.

\bibitem{FLZ}
{Vyacheslav Futorny, Zheng Li, and Jian Zhang}, {\it Irreducibility of the tensor product of Yangian\(\mathrm {Y}(\mathfrak {gl} _ {m| n})\) evaluation modules}, arXiv preprint arXiv:2607.24334 (2026).


\bibitem{FSZ}
{Vyacheslav Futorny, Vera Serganova, and Jian Zhang},
{\it Gelfand-Tsetlin modules for gl(m|n)},
Mathematical Research Letters 28.5 (2021): 1379-1418.


\bibitem{Gow1}
{L. Gow},
{\it  On the {Y}angian {$Y(\gl_{m|n})$} and its quantum {B}erezinian},
Czechoslovak J. Phys. {\bf 55} (2005), no. 11, 1415--1420.

\bibitem{Gow2}
{L. Gow},
{\it  Gauss decomposition of the {Y}angian {$Y(\gl_{m|n})$}},
Comm. Math. Phys. {\bf 276} (2007), no. 3, 799--825.

\bibitem{Lu}
{K. Lu},
{\it Gelfand-Tsetlin bases of representations for super Yangian and quantum affine superalgebra},
Lett. Math. Phys. {\bf 111} (2021), no.~6, Paper No. 145, 30 pp.

\bibitem{LM}
{K. Lu and E.~E. Mukhin},
{\it Jacobi-Trudi identity and Drinfeld functor for super Yangian},
Int. Math. Res. Not. IMRN {\bf 2021}, no.~21, 16751--16810.

\bibitem{Mol}
{A.~I. Molev},
Gelfand-Tsetlin basis for representations of Yangians, Lett. Math. Phys. {\bf 30} (1994), no.~1, 53--60.

\bibitem{Mol1}
{A.~Molev},
{\em Yangians for classical Lie algebras}.
Mathematical Surveys and Monographs
143. Amer. Math. Soc., Providence, RI, 2007.

\bibitem{Mol2}
{A. Molev},
{\it  Combinatorial bases for covariant representations of the {L}ie superalgebra {$\gl_{m|n}$}},
Bull. Inst. Math. Acad. Sin. (N.S.) {\bf 6} (2011), no. 4, 415--462.

\bibitem{MNO}
{A.~I. Molev, M. Nazarov and G.~I. Ol'shanski\u i},
{\it Yangians and classical Lie algebras},
Russian Math. Surveys {\bf 51} (1996), no.~2, 205--282; translated from Uspekhi Mat. Nauk {\bf 51} (1996), no.~2(308), 27--104.

\bibitem{Naz1}
{M. Nazarov},
{\it  Quantum {B}erezinian and the classical {C}apelli identity},
Lett. Math. Phys. {\bf 21} (1991), no.2, 123--131.

\bibitem{Naz2}
{M. Nazarov},
{\it Yangian of the general linear Lie superalgebra},
SIGMA Symmetry Integrability Geom. Methods Appl. {\bf 16} (2020), Paper No. 112, 24 pages.

\bibitem{NT1}
{M. Nazarov and V.~O. Tarasov},
{\it Yangians and Gel'fand-Zetlin bases},
Publ. Res. Inst. Math. Sci. {\bf 30} (1994), no.~3, 459--478.

\bibitem{NT2}
{M. Nazarov and V. Tarasov},
{\it  Representations of {Y}angians with {G}elfand-{Z}etlin bases},
J. Reine Angew. Math. {\bf 496} (1998), 181--212.

\bibitem{Palev1989a}
T.~D. Palev,
{\it Irreducible finite-dimensional representations of the Lie superalgebra {$\mathfrak{gl}(n/1)$} in a Gel'fand-Zetlin basis},
J. Math. Phys. {\bf 30} (1989), no.~7, 1433--1442.

\bibitem{Palev1989b}
T.~D. Palev,
{\it Essentially generic representations of the Lie superalgebras $\mathfrak{gl}(n|m)$ in the Gel'fand-Tsetlin basis},
Funct. Anal. Appl.
{\bf 23} (1989), no.~2, 141--142; translated from Funktsional. Anal. i Prilozhen.{\bf 23} (1989), no.~2, 69--70.

\bibitem{Ser}
{A.~N. Sergeev},
{\it Tensor algebra of the identity representation as a module over the Lie superalgebras ${\rm Gl}(n,\,m)$\ and $Q(n)$},
Mat. Sb. (N.S.) {\bf 123(165)} (1984), no.~3, 422--430.

\bibitem{SV}
{N. Stoilova and J. Van der Jeugt},
{\it  Gel'fand-{Z}etlin basis and {C}lebsch-{G}ordan coefficients for covariant representations of the {L}ie superalgebra {$\gl(m|n)$}},
J. Math. Phys. {\bf 51} (2010), no. 9, 093523, 15.

\bibitem{Zhang1}
{R. Zhang},
{\it Representations of super Yangian},
J. Math. Phys. {\bf 36} (1995), no.~7, 3854--3865.

\bibitem{Zhang2}
{R. Zhang},
{\it  The {${\rm gl}(M|N)$} super {Y}angian and its finite-dimensional representations},
Lett. Math. Phys. {\bf 37} (1996), no. 4, 419--434.

\end{thebibliography}
\end{document}